\input amstex

\define\ab{\text{\hglue.2mm}}
\define\eps{\varepsilon}
\loadbold \loadeufm \loadmsbm

\newsymbol\gtrsim 1326
\newsymbol\lesssim 132E

\define\noi{\noindent{\hglue2.5mm}}
\define\m{\medskip}
\define\s{\smallskip}
\define\te{\text{ }}
\define\hm{\text{\hglue .5mm}}

\define\col{\!\ab:\ab}

\define\PP{\text{\rm P}}
\define\EE{\text{\rm E}}

\define\num#1{{$\boxed{\bold{#1}}$\ab}}

\magnification=1200

%\rightline{\trr Version of \ab6/25/2020}

\m

\centerline{\bf Diffusion approximation for noise-induced evolution}

\centerline{\bf of first integrals in multifrequency systems}

\bigskip

\centerline{M.I.\,Freidlin, \qquad A.D.\ab Wentzell}

\bigskip

\centerline{\bf Abstract}

\m

We consider fast oscillating random perturbations of dynamical 
systems in regions where one can introduce action-angle-type 
coordinates. In an appropriate time scale, the evolution of
first integrals, under the assumption that the set of resonance
tori is small enough, is approximated by a diffusion process.
If action-angle coordinates can be introduced only piece-wise,
the limiting diffusion process should be considered on an 
open-book space. Such a process can be described by differential
operators, one in each page, supplemented by some gluing 
conditions at the binding of the open book.

\m

\noi{\bf Key words}: averaging principle, diffusion approximation,
random perturbations, resonance tori.

\m

\noi{\bf AMS subject classification}: 37J40, 60HXX, 58G32.

\bigskip

\noindent{\bf 0. Introduction.}

\m

Let $\,\boldsymbol\xi^\eps(t)$, $\ab\eps > 0$, $\,t \geq 0$,
\ab be the solution of the $N$-dimensional differential equation
$$
    \dot{\boldsymbol \xi}^\eps(t)
    = \boldkey b\bigl(\boldsymbol\xi^\eps(t),\, \zeta_{t/\eps}\bigr),
    \quad \boldsymbol\xi^\eps(0) = \boldkey x_0,
\eqno(0.1)
$$
where $\,\zeta_s\ab$ is a stationary stochastic process on an
arbitrary space $\Cal Z$ ($0 <\eps << 1$).
If the process $\,\zeta_s\ab$ has good enough mixing properties
(and under some conditions on $\,\boldkey b(\boldkey x,\, z)$),
the stochastic process $\,\boldsymbol\xi^\eps(t)\ab$ will converge
in probability, uniformly on every finite time interval, to
the solution of the averaged differential equation
$$
    \dot{\boldkey X}(t) = \tilde{\boldkey b}\bigl(\boldkey X(t)\bigr),
    \quad \boldkey X(0) = \boldkey x_0,
\eqno(0.2)
$$
where
$$
    \tilde{\boldkey b}(\boldkey x)
    = \EE\, \boldkey b(\boldkey x,\, \zeta_s)
    = \int_{\Cal Z} \boldkey b(\boldkey x,\, z)\ \mu(dz),
\eqno(0.3)
$$
$\,\mu\,$ being the distribution of $\,\zeta_s$ -- not depending
on $\,s$ (see, for instance, [7], Theorem 7.2.1 and
the references there). This means that $\,\boldsymbol\xi^\eps(t)\ab$
can be considered as the result of a small random perturbation
of $\boldkey X(t)$; this is a result of the Law-of-Large-Numbers
type. A similar result holds if $\,\zeta_s\ab$ is not a stationary
process, but a stochastic process with good mixing properties
whose distribution at time $\,s\,$ approaches a limiting
distribution $\,\mu\,$ as~$\,s \to \infty$.

Let us denote $\,\boldsymbol\alpha(\boldkey x,\, z) =
\boldkey b(\boldkey x,\, z) - \tilde{\boldkey b}(\boldkey x,\, z)$;
\ab the integral of $\,\boldsymbol\alpha(\boldkey x,\, z)\ab$
is equal to $\bold0$. This is the perturbation function.

\s

One may be interested in what the behavior of $\,\boldsymbol
\xi^\eps(t)\ab$ is on time intervals whose length goes to~$\infty$
as $\,\eps \to 0\ab$: \ab the behavior of $\,\boldsymbol\xi^\eps
\bigl(t(\eps)\bigr)$, $\ab t(\eps) \to \infty$ ($\eps \to 0$).
Of course it depends on the rate at which $\,t(\eps)\ab$ goes
to $\infty$. It turns out that, typically, significant deviations
of $\,\boldsymbol\xi^\eps(t)\ab$ from $\boldkey X(t)$ occur
on time intervals of length of order of $\,\eps^{- 1}$; \ab
one can be interested in the limiting behavior of the stochastic
process~$\boldkey X^\eps(t) = \boldsymbol\xi^\eps(t/\eps)$,
obtained from our original process by a change of time parameter.
This process is the solution of the equation
$$
    \dot{\boldkey X}^\eps(t)
    = \eps^{- 1} \cdot \boldkey b\bigl(\boldkey X^\eps(t),\,
    \zeta_{t/\eps^2}\bigr), \quad \boldkey X^\eps(0) = \boldkey x_0.
\eqno(0.4)
$$

If the system (0.2) has a first integral $H(\boldkey x)$
(so that $H\bigl(\boldkey X(t)\bigr) \equiv H(\boldkey x_0)$),
then $H\bigl(\boldsymbol\xi^\eps(t)\bigr)$ converges in probability
to the constant $H(\boldkey x_0)$ in every finite interval of~
$\,t$'s. \ab There may\linebreak be $\,n\,$ independent first integrals
$\ab H_1(\boldkey x)$, ..., $H_n(\boldkey x)$; \ab let us introduce
the vector-valued function $\boldkey H(\boldkey x) = \bigl(
H_1(\boldkey x),\ab...,\,H_n(\boldkey x)\bigr)$.

\s

For the averaged system (0.2) having first integrals the
perturbed system (0.!) (or (0.4)) may have first
integrals or it may not.

\s

In the particular case of $\,n = N$, \ab all coordinates $\,x_1$,
..., $x_N\ab$ of $\,\boldkey x\,$ being first integrals (which
means that the system (0.2) is just $\dot{\boldkey X}(t) =
\bold0$), the results on the limiting behavior of
$\boldkey X ^\eps(t)$ were obtained in [11], [9],
[3];
these results were not about convergence in probability, but
rather (as it should be) about {\it convergence in distribution\/}.
i.\.a., weak convergence of distributions of the trajectories
$\boldkey X^\eps(\bullet)$ in the space $\bold C[0,\, T]$ of
continuous functions for every $T \in [0,\, \infty)$.
In the case of $\,n < N\ab$ it
may happen that all $\,p$-dimensional ($p = N - n$) level surfaces
$\{\boldkey x\col \boldkey H(\boldkey x) = \boldkey h\}$ within
a region $G$ in our $N$-dimensional space are diffeomorphic
to each other. Then we can introduce new coordinates in $G$:
an\linebreak $\,n$-dimensional \ab coordinate $\,\boldkey h = \boldkey H(
\boldkey x)$, and a  $\,p$-dimensional $\,\boldkey y \in \Cal
Y$, \ab where $\Cal Y$ is the manifold diffeomorphic to all
level surfaces $\{\boldkey x\col \boldkey H(\boldkey x)
= \boldkey h\}$ in our region. Let us denote $\boldkey H^\eps(t)$
the $\,\boldkey h$-coordinate \ab of $\boldkey X^\eps(t)$, and
$\boldkey Y  ^\eps(t)$ its $\,\boldkey y$-coordinate. In these
new coordinates (0.4) takes the form
$$
\left\{\aligned
    \dot{\boldkey H}^\eps(t) &    = \eps^{- 1} \cdot \boldsymbol
    \beta_{\boldkey h}\bigl(\boldkey H^\eps(t),\, \boldkey Y^\eps(t),\,
    \zeta_{t/\eps^2}\bigr), \\
    \dot{\boldkey Y}^\eps(t) &    = \eps^{- 1} \cdot\boldsymbol
    \beta_{\boldkey y}\bigl(\boldkey H^\eps(t),\, \boldkey Y^\eps(t),\,
    \zeta_{t/\eps^2}\bigr),
\endaligned\right.
\eqno(0.5)
$$
where $\,\boldsymbol\beta_{\boldkey y}\ab$ takes values in the
tangent bundle
of the manifold $\Cal Y$. Of the corresponding averaged vector
fields
$$
    \tilde{\boldsymbol\beta}_{\boldkey h}(\boldkey x)
    = \int_{\Cal Z} \boldsymbol\beta_{\boldkey h}(\boldkey x,\, z)\ \mu(dz), \quad
    \tilde{\boldsymbol\beta}_{\boldkey y}(\boldkey x)
    = \int_{\Cal Z} \boldsymbol\beta_{\boldkey y}(\boldkey x,\, z)\ \mu(dz)
\eqno(0.6)
$$
the first is equal to $0$, while the second one is, generally,
non-zero. So we can call the $\,\boldkey h$-coordinates \ab``slow''
coordinates, and $\,\boldkey y$-coordinates ``fast'' ones.

For $\,n = N - 1$, $p = 1\ab$ the level manifolds $\{\boldkey x
\col \boldkey H(\boldkey x) = \boldkey h\}$, if they are compact,
are diffeomorphic to a {\it circle\/}. If
$\,\beta_y(\boldkey x) \neq 0$, one can change the coordinates
so that the new $\,\tilde\beta_y(\boldkey x)\ab$ depends
only on $\,\boldkey h$: $\ab\tilde\beta_y(\boldkey x) =
\omega(\boldkey h) \neq 0\ab$ being the {\it frequency\/}
of the rotation of $\boldkey X(t)$ on the level circle $\{\boldkey x
\col \boldkey H(\boldkey x) = \boldkey h\}$ (in this case the
$\,y$-coordinate \ab can be considered as an angle and denoted
with letter $\,\varphi$, \ab the corresponding component of
$\boldkey X(t)$, $\boldkey X^\eps(t)$ being denoted $\Phi(t)$,
$\Phi^\eps(t)$). \ab In the case of $\,n = p = 1\ab$ results
about convergence in distributions of the ``slow'' component of the
process $\boldkey X^\eps(t)$ were obtained in [4]. One should
also mention [5], where this type of questions were studied

\s

Our goal is to study the problem for $\,p > 1$. We restrict
ourselves to the case in which all level surfaces are $\,p$-dimensional
\ab unit tori, and the equation (0.2) is
$$
\left\{\aligned
    \dot{\boldkey H}(t) &= \bold 0,       \\
    \dot{\boldsymbol\Phi}(t) &= \boldsymbol\omega\bigl(
            \boldkey H(t)\bigr)
\endaligned\right.
\eqno(0.7)
$$
($\boldsymbol\omega(\boldkey h) = \bigl(\omega_1(\boldkey h),\ab
...,\, \omega_p(\boldkey h)\bigr)$ is the {\it vector\/} of
frequencies).

If, for a fixed $\,\boldkey h$, \ab the frequencies $\,\omega_j
(\boldkey h)\ab$ are rationally independent, there is only
one invariant measure for the system (0.7) on the torus
$\{\boldkey x\col \boldkey H(\boldkey x) = \boldkey h\}$ with
the total value equal to $1$, namely (in the $\,\boldsymbol
\varphi $-coordinates), the $\,p$-dimensional Lebesgue measure;
and there is some mixing for this system: the
time average of a continuous function of~$\,\boldsymbol\Phi
(t)\ab$ over a growing time interval converges to the average
of this function over the torus. Something similar should be
true for the stochastic process $\boldsymbol\Phi^\eps(t)$, only
time intervals for~$\boldsymbol\Phi(t)$ of length going to
$\infty$ may correspond to time intervals for $\boldsymbol\Phi
^\eps(t)$ of infinitely small lengths.

If $\,\omega_j(\boldkey h)\ab$ are rationally dependent:
$$
    \sum_{\,j = 1}^p  k_j \cdot \omega_j(\boldkey h) = 0,
\eqno(0.8)
$$
$k_j\ab$ being integers, $\ab\boldkey k = (k_1,\ab...,\, k_p)
\neq \bold 0$\ab, there are infinitely many invariant measures
with total value $1$; so averaging over the Lebesgue measure
seems to be for nothing here. If the equality (0.8) is
satisfied for all $\,\boldkey h$, \ab we can add to $\ab H_1
(\boldkey x)$, ..., $H_p(\boldkey x)\ab$ another first integral:
$$
    H_{p + 1}(\boldkey x)
    = \sum_{j = 1} ^p k_j \cdot \varphi_j
\eqno(0.9)
$$
(the values of the coordinates $\,\varphi\ab$ are taken in the
interval $[\ab 0,\, 1)$;
if at least one of $\,k_j\ab$ were not an integer, the function
(0.9) would be not smooth).

Normally, tori with rationally independent and rationally dependent
frequencies $\,\omega_j(\boldkey h)\ab$ alternate; and the
natural assumption under which one should try to obtain results
in our problem is that the set of $\,\boldkey h\,$ for which
(0.7) holds is small in some sense; the most natural assumption
would be that the lebesgue measure of the set of $\,\boldkey h$'s
\ab for which the frequencies are rationally dependent is equal
to $0$.

\s

We are considering unperturbed systems of the form (0.7)
not only because of the convenience of using only one local
coordinate system (the tangent bundle to $\Bbb T^p$ being
identified with the space $\Bbb R^p$) but because of some deeper
reasons: level surfaces diffeomorphic to tori appear naturally
for some Hamiltonian systems.

\s

Let $N = 2\ab n\ab$, $\ab \boldkey x = (\boldkey p,\, \boldkey q)
= (p_1,\ab...,\, p_n; \, q_1,\ab...,\, q_n)$ \ab Let $H(\boldkey
x)$ be a smooth enough function; $\overline\nabla H(\boldkey x) =
\bigl(- \nabla_{\boldkey q}\ab H(\boldkey p,\, \boldkey q),\,
\nabla_{\boldkey p}\, H(\boldkey p,\, \boldkey q)\bigr)$. Consider
the system
$$
    \dot{\boldkey X}(t) = \overline\nabla H\bigl(\boldkey X(t)\bigr),
    \quad \boldkey X(0) = \boldkey x_0\ab.
\eqno(0.!0)
$$
Assume that the Hamiltonian system (0.10) is completely integrable,
which means that it has $\,n\,$ smooth first integrals $H_1(\boldkey x)
= H(\boldkey x)$, $H_2(\boldkey x)$, ..., $H_n(\boldkey x)$ such that
$\nabla H_i(\boldkey x) \cdot \overline\nabla H_j(\boldkey x)$
$= 0$ for all $\,i$, $j$; \ab the level sets $\{\boldkey x\col
H_1(\boldkey x) = h_1,\ab...,\, H_n(\boldkey x) = h_n\}$ in
some region $G$ are connected compact smooth manifolds; and
the gradients $\nabla H_i(\boldkey x)$ are linearly independent
for every $\,\boldkey x \in G$. Then (see [1], Section
49) the level sets are diffeomorphic to the\linebreak
$\,n$-dimensional torus $\Bbb T ^n$, and one can introduce
in the region $G$ coordinates $\ab\boldkey H$, $\boldsymbol\varphi
\ab$ ({\!\it action-angle coordinates\/}) such that the system
(0.10) takes the form (0.7) with  $\,p = n$.

\s

We are going to study our problem under some simplifying assumptions.
First, that the space in which the slow coordinate $\,\boldkey h\,$
changes is the whole Euclidean space $\Bbb R^n$ (so that $\,
\boldkey x\,$ changes in $\Bbb R^n \times \Bbb T^p$). If we are
thinking of applying our results, say, to the perturbations
of completely integrable Hamiltonian systems, we have to consider
the perturbed system only up to the random time at which it
leaves the region $G$. For example, the action-angle coordinates
don't work in neighborhoods of critical points of the function
$H(\boldkey x)$. So to apply our results we have to be
prepared to the action variables running over some region that
is smaller than the whole Euclidean space. This is not a very
serious restriction: if we have system (0.1) (or (0.4))
only in a bounded
region, we can extend its coefficients to the whole space, apply
our results, and then consider the limiting process only up to
the time at which it leaves the region in question. Considering
processes only in some region can be done in the same way as
in the present paper; but with some bother, e.\, g., we have to consider
integrals with random times as their upper limit,
etc.

Of course it's a more difficult problem to study the perturbed
system in regions that do contain critical points and the like
(see Subsection 3.2);
our present paper can be considered as the first step in this direction.

\s

The second simplifying restriction is that the ``driving''
stochastic process $\,\zeta_s\ab$ is a
non-degenerate finite-dimensional diffusion process on a compact
manifold~$\Cal Z$. It turns out in this case that the trajectories
of the stochastic process $\boldkey X^\eps(t)$ are close to
those of a stochastic process $\tilde{\boldkey X}^\eps(t)$
that is expressed by means of integrals and stochastic integrals
(with the integrands in the integrals for $\,h$-components
being $O(1)$)\ab:
not precisely a diffusion process, but rather a component of
one of a higher dimension. This allows us to use stochastic-equations
technique all over our reasoning~-- instead of a combination
of this technique with that used for sums of dependent random
variables. After this, we don't need to turn to the mixing
properties of the ``driving'' process (they are
anyway very good: the dependence between its values that are
separated by a large time interval decreases exponentially
as this interval grows).

The problem of asymptotic behavior of a diffusion process with
large drift was considered in [7],
Chapters 8, 9; and our present problem is considered
in a similar way.

Further simplifying assumptions: that the manifold  $\Cal Z$
is an $\,m$-dimensional unit torus~$\Bbb T^m$, which allows
us not to think about changing from one map on the manifold
to another\ab; and that the process $\,\zeta_s\ab$ is the standard
Wiener process on $\Bbb T^m$.

\s

So we have two tori: $\Bbb T^p$ and $\Bbb T^m$,
but their status is different: the first one is because of the
action-angle-type coordinates (without which it would be a completely
different problem); the other because of our simplifying
assumptions.

\bigskip

\noindent{\bf 1. Averaging over $\Bbb T^m$. Formulation of the
main result.}

\m

In the Introduction we outlined the problem we are going to
consider; let us repeat it here, with precise formulations.

\s

Let $\,\boldkey b(\boldkey x,\, \boldkey w)$, $\ab\boldkey x
\in \Bbb R^n \times \Bbb T^p$, $\boldkey w \in \Bbb T^m$, \ab
be a function with values in $\Bbb R^N$ ($N = n + p$), bounded
and continuous in $(\boldkey x,\,
\boldkey w)$ and H\"older-continuous in $\,\boldkey w$, \ab
uniformly in $\,\boldkey x$, together with its derivatives
$\dfrac{\partial \boldkey b}{\partial x_i}$, $\dfrac{\partial^2
\boldkey b}{\partial x_i\ab \partial x_j}$,  $\dfrac{\partial^3
\boldkey b}{\partial x_i\ab \partial x_j\ab \partial x_k}$.

Let ${\boldkey w}(t)$ be a standard Wiener process on $\Bbb
T^m$; let $\boldkey W^\eps(t) = {\boldkey w}(t/\eps^2)$.
Let $\boldkey X^\eps(t)$ be the stochastic process defined
as the solution of the differential equation
$$
    \dot{\boldkey X}^\eps(t)
    = \eps^{- 1}\ab \boldkey b\bigl(\boldkey X^\eps(t),\,
            \boldkey W^\eps(t)\bigr)
\eqno(1.1)
$$
with an initial condition $\boldkey X^\eps(0) = \boldkey x_0$
not depending on $\,\eps$.
The process $\bigl(\boldkey X^\eps(t),\, \boldkey W^\eps(t)\bigr)$
is the solution of the system of stochastic equations
$$
    d\ab \boldkey X^\eps(t) = \eps^{- 1}\ab
    \boldkey b\bigl(\boldkey X^\eps(t),\, \boldkey W^\eps(t)\bigr)\, dt\ab,
    \quad d\ab \boldkey W^\eps(t) = \eps^{- 1} \ab
    d\ab \boldkey W(t),
\eqno(1.2)
$$
where $\boldkey W(t)$ is another standard Wiener process. The
first equation here means that
$$
    \boldkey X^\eps(t) = \boldkey x_0
    + \eps^{-1} \int_0^t \boldkey b\bigl(\boldkey X^\eps(s),\,
        \boldkey W^\eps(s)\bigr)\ ds.
\eqno(1.3)
$$

If a H\"older-continuous function $G(\boldkey w)$ is such
that $\dsize\int_{\Bbb T^m} G(\boldkey w)\ d\boldkey w = 0$,
there exists a solution $U(\boldkey w)\ab$ of the equation $\tfrac12\ab
\Delta\ab U(\boldkey w) = -\, G(\boldkey w)$, $\ab\boldkey w \in \Bbb T^m$,
and it is unique up to an additive constant.

\m

{\bf Lemma 1.1.} {\it Let\/ $G(\boldkey x,\, \boldkey w),$
$\ab \boldkey x \in \Bbb R^N,$ $\boldkey w \in \Bbb T^m,$ \ab
be a function that is$,$ together with its\/ $\,x$-derivatives
\ab up to order\/ $K,$ bounded$,$ continuous in\/ $(\boldkey x,\,
\boldkey w)$ and H\"older continuous in\/ $\,\boldkey w,$ \ab
uniformly in\/ $\,\boldkey x;$ \ab let\/ $\dsize\int_{\Bbb T^m}
G(\boldkey x,\, \boldkey w)\ d\boldkey w = 0$ for every
$\,\boldkey x \in \Bbb R^N.$

\s

Then there exists a solution\/ $U(\boldkey x,\, \boldkey w)$ of the
equation
$$
    \tfrac12\ab \Delta_{\boldkey w}\ab U(\boldkey x,\, \boldkey w)
    = -\,G(\boldkey x,\, \boldkey w), \quad
    \boldkey x \in \Bbb R^N, \ \boldkey w \in \Bbb T^m,
\eqno(1.4)
$$
that is, together with its derivatives $\botsmash{\dfrac{\partial U}
{\partial  w_j}}$, $\botsmash{\dfrac{\partial^k U}{\partial
x_{i_1}...\,\partial x_{i_k}}}$, $\botsmash{\dfrac{\partial^{\ab k + 1} U}
{\partial w_j\ab\partial x_{i_1}...\,\partial x_{i_k}}}$, $\ab k \leq K$,
bounded and continuous}.

\m

The {\bf proof} can be based on the ``probabilistic'' representations:
$$
\aligned
    U(\boldkey x,\, \boldkey w)
    = &\int_0^\infty \EE_{\boldkey w}\, G\bigl(\boldkey x,\,
        \boldkey w(t)\bigr)\ dt, \quad
    \dfrac{\partial U}{\partial w_j}(\boldkey x,\, \boldkey w)
    = \int_0^\infty \dfrac{\partial}{\partial w_j}
     \EE_{\boldkey w}\, G\bigl(\boldkey x,\,
        \boldkey w(t)\bigr)\ dt,        \\
    \dfrac{\partial^k U}{\partial x_{i_1}...\,\partial x_{i_k}}
    &(\boldkey x,\, \boldkey w)
    = \int_0^\infty \EE_{\boldkey w}\, \dfrac{\partial^k G}{
    \partial x_{i_1}...\,\partial x_{i_k}} \bigl(\boldkey x,\,
        \boldkey w(t)\bigr)\ dt, \\
    &\dfrac{\partial^{\ab k + 1} U}{\partial w_j\ab
        \partial x_{i_1}...\,\partial x_{i_k}}
    (\boldkey x,\, \boldkey w)
    = \int_0^\infty \dfrac{\partial}{\partial w_j}
    \EE_{\boldkey w}\, \dfrac{\partial^k G}{
    \partial x_{i_1}...\,\partial x_{i_k}} \bigl(\boldkey x,\,
        \boldkey w(t)\bigr)\ dt,
\endaligned
\eqno(1.5)
$$
where $\ab\EE_{\boldkey w}$ denotes the expectation evaluated
under the assumption that the Wiener process starts from the
point $\,\boldkey w(0) = \boldkey w$.

\m

{\bf Lemma 1.2.} {\it Let a function\/ $\,g(\boldkey x,\, \boldkey w)\ab$
be bounded and continuous in $(\boldkey x,\, \boldkey w)$ and
H\"older-continuous in\/ $\,\boldkey w\,$ together with its
derivatives\/ $\botsmash{\dfrac{\partial g}{\partial x_i}}.$
Let us define the function
$$
    \tilde g(\boldkey x)
    = \botsmash{\int_{\Bbb T^m} g(\boldkey x,\, \boldkey w)\ d\boldkey w.}
\eqno(1.6)
$$

Then
$$
\aligned
    \eps^{- 1} \cdot\int_0^t & g(\boldkey X^\eps(s),\,
            \boldkey W^\eps(s)\bigr)\ ds
    - \eps^{-  1} \cdot \int_0^t \tilde g\bigl(
            \boldkey X^\eps(t)\bigr)\ ds \\
    &\quad= \int_0^t \nabla_{\boldkey w} U\bigl(
     \boldkey X^\eps(s),\, \boldkey W^\eps(s)\bigr)\ d\ab \boldkey W(s)\\
     &\qquad\qquad+ \int_0^t \nabla_{\boldkey x} U\bigl(
     \boldkey X^\eps(s),\, \boldkey W^\eps(s)\bigr) \cdot
     \boldkey b\bigl(\boldkey X^\eps(s),\, \boldkey W^\eps(s)\bigr)\ ds
    \ +\ O(\eps)\ab,
\endaligned
\eqno(1.7)
$$
where the\/ $O(\eps)\!\ab$ holds uniformly}.

\m

{\bf Proof.} The function $G(\boldkey x,\, \boldkey w) =
g(\boldkey x,\, \boldkey w) - \tilde g(\boldkey x)\ab$ satisfies
the conditions of Lemma~1.1 with $K = 1$. Formula (1.7) (with
$O(\eps) = \eps \cdot
\bigl[U\bigl(\boldkey X^\eps(0),\, \boldkey W^\eps(0)\bigr)
- U\bigl(\boldkey X^\eps(t),\, \boldkey W^\eps(t)\bigr)\bigr]$)
is obtained by applying the It\^o formula to $\,\eps \cdot U\bigl(
\boldkey X^\eps(t),\, \boldkey W ^\eps(t)\bigr)$.

\m

{\bf Lemma 1.3.} {\it Under the conditions of\/ {\rm Lemma 1.2}
we have\ab}:
$$
    \EE\ab\Bigl|\int_0^t g(\boldkey X^\eps(s),\,
            \boldkey W^\eps(s)\bigr)\ ds
    - \int_0^t \tilde g\bigl(
            \boldkey X^\eps(s)\bigr)\ ds\ab\Bigr|
    = O(\eps),
\eqno(1.8)
$$
{\it with\/ $O(\eps)\!$ being uniform in\/ $\,t\,$ changing in
every finite interval\ab}.

\s

{\bf Proof}: Multiply (1.7) by $\,\eps\,$ and take the expectation
of the absolute value of both sides.

\s

Sometimes we'll need the versions of formulas like (1.7),
(1.8) with integrals from $0$ to $\infty$; which is possible
if we multiply the integrands by $\,e ^{- \lambda s}$, $\ab
\lambda > 0$.

\s

{\bf Lemma 1.3$^{\ab\boldsymbol\prime}$.} {\it Under the conditions
of\/ {\rm Lemma 1.2} for every\/} $\,\lambda > 0\ab$
$$
    \EE\ab\Bigl|\int_0^\infty e^{- \lambda s} \hm
    \bigl[g(\boldkey X^\eps(s),\, \boldkey W^\eps(s)\bigr)
    - \tilde g\bigl(\boldkey X^\eps(s)\bigr)\bigr]\ ds\ab\Bigr|
    = O(\eps).
\eqno(1.9)
$$

The {\bf proof} is obtained by applying the It\^o formula to
$\,\,\eps^2 \cdot e^{- \lambda t} \cdot  U\bigl(\boldkey X^\eps(t),
\, \boldkey W^\eps(t)\bigr)$, taking the expectation of the
absolute value, and letting $\,t \to \infty$.

\m

Now let the function $\,\tilde{\boldkey b}(\boldkey x)\ab$
be defined by
$$
    \tilde{\boldkey b}(\boldkey x)
    = \int_{\Bbb T^m} \boldkey b(\boldkey x,\, \boldkey w)\
        d\ab\boldkey w
\eqno(1.10)
$$
(which replaces formula (0.3)); $\ab\boldsymbol\alpha(\boldkey x,\,
\boldkey w) = \boldkey b(\boldkey x,\, \boldkey w) - \tilde{
\boldkey b}(\boldkey c)$. Let $\,\boldkey u(\boldkey x,\,
\boldkey w)\ab$ be a bounded solution of the equation
$$
    \tfrac12\ab \Delta_{\boldkey w} \boldkey u(\boldkey x,\, \boldkey w)
    = -\, \boldsymbol\alpha(\boldkey x,\, \boldkey w), \quad
    \boldkey w \in \Bbb T^m.
\eqno(1.11)
$$
Applying Lemma 1.2 to the {\it vector\/} function $\,\boldkey u(
\boldkey x,\, \boldkey w)$, we get:
$$
\aligned
    \eps^{- 1}\int_0^t & \boldkey b(\boldkey X^\eps(s),\,
            \boldkey W^\eps(s)\bigr)\ ds
    - \eps^{- 1} \int_0^t \tilde{\boldkey b}\bigl(
                \boldkey X^\eps(t)\bigr)\ ds \\
    &\quad= \int_0^t \nabla_{\boldkey w} \boldkey u\bigl(
     \boldkey X^\eps(s),\, \boldkey W^\eps(s)\bigr)\ d\ab \boldkey W(s)\\
     &\qquad\qquad+ \int_0^t \nabla_{\boldkey x} \boldkey u\bigl(
     \boldkey X^\eps(s),\, \boldkey W^\eps(s)\bigr) \cdot
     \boldkey b\bigl(\boldkey X^\eps(s),\, \boldkey W^\eps(s)\bigr)\ ds
    \ +\ O(\eps)
\endaligned
\eqno(1.12)
$$
($\nabla_{\boldkey w} \boldkey u\,$ and $\ab\nabla_{\boldkey x}
\boldkey u\,$ are matrices of sizes $\ab N \times m\,$ and
$N \times N$, correspondingly). From this and formula (1.3)
we get:
$$
\aligned
    \boldkey X^\eps&(t) = \boldkey x_0
    + \topsmash{\int_0^t} \nabla_{\boldkey w}\ab
      \boldkey u\bigl(\boldkey X^\eps(s),\,
        \boldkey W^\eps(s)\bigr)\ d\ab\boldkey W(s)     \\
    &+ \int_0^t \bigl[\ab \eps^{- 1} \cdot \tilde{\boldkey b}\ab
     \bigl(\boldkey X^\eps(s)\bigr)
    + \nabla_{\boldkey x}\ab
      \boldkey u\bigl(\boldkey X^\eps(s),\,
      \boldkey W^\eps(s)\bigr)\cdot \boldkey b\ab\bigl(
       \boldkey X^\eps(s),\,
        \boldkey W^\eps(s)\bigr)\bigr]\ ds + O(\eps).
\endaligned
\eqno(1.13)
$$
So $\boldkey X^\eps(t)$ is close~-- in the sense of uniform
closeness of trajectories~-- to the stochastic process
$$
\aligned
    \tilde{\boldkey X}^\eps(&t) = \boldkey x_0
    + \topsmash{\int_0^t} \nabla_{\boldkey w}\ab
      \boldkey u\bigl(\boldkey X^\eps(s),\,
        \boldkey W^\eps(s)\bigr)\ d\ab\boldkey W(s)     \\
    &+ \int_0^t \bigl[\ab\eps^{- 1} \cdot \tilde{\boldkey b}\ab
     \bigl(\boldkey X^\eps(s)\bigr)
    + \nabla_{\boldkey x}\ab
      \boldkey u\bigl(\boldkey X^\eps(s),\,
      \boldkey W^\eps(s)\bigr)\cdot \boldkey b\ab\bigl(
       \boldkey X^\eps(s),\,
        \boldkey W^\eps(s)\bigr)\bigr]\ ds\ab.
\endaligned
\eqno(1.14)
$$
This process is an $N$-dimensional
component of a diffusion process of a higher dimension; one
component of a multidimensional diffusion process need not be
a diffusion process, so we are not close to establishing closeness
of $\boldkey X^\eps(t)$ to one.

\m

Now let us look what happens if $(\boldkey h,\, \boldsymbol\varphi)$
are action-angle-type coordinates: $\ab\boldkey x = (\boldkey h,\,
\boldsymbol\varphi)$, $\boldkey h \in \Bbb R^n$, $\boldsymbol\varphi
\in \Bbb T^p$, $p = N - n$, $\ab \tilde{\boldkey b}(\boldkey x) = \tilde{
\boldkey b}(\boldkey h,\, \boldsymbol\varphi) = \bigl(\bold 0,\,
\boldsymbol\omega(\boldkey h)\bigr)$, $\ab \boldkey X^\eps(t) = \bigl(
\boldkey H ^\eps(t),\, \boldsymbol\Phi ^\eps(t)\bigr)$.

\m

The equations (1.13) are separated into those for slow and
fast components:
$$
    \boldkey H^\eps(t) = \tilde{\boldkey H}^\eps(t) + O(\eps),
    \qquad \boldsymbol\Phi^\eps(t)
    = \tilde{\boldsymbol\Phi}^\eps(t) + O(\eps),
\eqno(1.15)
$$
$$
\aligned
    \tilde{\boldkey H}^\eps(&t) = \boldkey h_0
    + \topsmash{\int_0^t} \nabla_{\boldkey w}\,
      \boldkey u_{\boldkey h} \bigl(\boldkey H^\eps(s),\,
      \boldsymbol\Phi^\eps(s),\,
        \boldkey W^\eps(s)\bigr)\ d\ab\boldkey W(s)     \\
    &+ \int_0^t \nabla_{\boldkey x}\,
      \boldkey u_{\boldkey h}\bigl(\boldkey H^\eps(s),\,
      \boldsymbol\Phi^\eps(s),\,
      \boldkey W^\eps(s)\bigr)\cdot \boldkey b\ab\bigl(
      \boldkey H^\eps(s),\, \boldkey Z^\eps(s),\,
        \boldkey W^\eps(s)\bigr)\ ds\,,
\endaligned
\eqno(1.16)
$$
$$
\aligned
    \tilde{\boldsymbol\Phi}^\eps(t) & = \boldsymbol\varphi_0
    + \topsmash{\int_0^t} \nabla_{\boldkey w}\,
      \boldkey u_{\boldsymbol\varphi}\bigl(\boldkey H^\eps(s),\,
      \boldsymbol\Phi^\eps(s),\,
        \boldkey W^\eps(s)\bigr)\ d\ab\boldkey W(s)
        + \eps^{- 1} \int_0^t \boldsymbol\omega\ab \bigl(
            \boldkey H^\eps(s)\bigr)\ ds    \\
    &+\! \int_0^t \nabla_{\boldkey x}\,
      \boldkey u_{\boldsymbol\varphi}\bigl(\boldkey H^\eps(s),\,
      \boldsymbol\Phi^\eps(s),\,
        \boldkey W^\eps(s)\bigr)\cdot \boldkey b\ab\bigl(
       \boldkey H^\eps(s),\, \boldsymbol\Phi^\eps(s),\,
        \boldkey W^\eps(s)\bigr)\ ds\ab,
\endaligned
\eqno(1.17)
$$
where the subscripts $\,_{\boldkey h}$\ab, $_{\boldsymbol\varphi}\,$
mean the $\,\boldkey h$-, $\boldsymbol\varphi$-components
\ab of the vector-valued function $\,\boldkey u\ab$. \ab We'll
denote the separate one-dimensional components of
$\,\boldkey u_{\boldkey h}\ab$ as $\,u_i$\ab, $\ab i = 1$, ..., $n$.

\s

Of course one cannot expect that the process
$\boldsymbol\Phi^\eps(t)$ should converge in distribution
as $\,\eps \to 0\ab$ to anything\ab; but we can
expect that $\boldkey H^\eps(t)$ and $\tilde{\boldkey H}
^\eps(t)$ converge. Because of (1.15) their weak limits
must be the same.

The standard way to establish weak convergence of function-space
distributions is to first establish {\it tightness\/} (weak
pre-compactness) of the family of distributions.

\s

If the distribution of one of $\boldkey H^\eps(\bullet)$,
$\tilde{\boldkey H}^\eps(\bullet)$ converges weakly as $\,\eps \to 0$, \ab while
the family of distributions of the other is tight, the other
distribution also converges weakly, and to the same limit.

\m

{\bf Lemma 1.4.} {\it Let\/ $\,\xi^\eps(t),$ $t \in [0,\, T],$
be a family of stochastic processes with values in a complete
metric space with distance\/ $\,\rho(\ \ ,\ \ ).$ Let the following
inequality hold for all\/ $\,t,$ $s \in [0, T]$ \ab and\/ $\,\eps\!:$
$$
    \EE\,\rho\bigl(\xi^\eps(t),\, \xi^\eps(s)\bigr)^\beta
    \leq |t - s|^{1 + \alpha},
\eqno(1.18)
$$
where $\,\alpha,$ $\beta\,$ are positive constants\ab$;$ and let
the one-dimensional distributions of\/ $\,\xi^\eps(0)\ab$ form
a tight family$.$

Then the family of distributions of the trajectories $\,\xi^\eps(
\bullet)\ab$ in the space\/ $\Bbb C[0,\, T]$ of continuous
functions is tight}.

\s

This is an adaptation of Kolmogorov's theorem about the continuous
modification: see [8], Theorem 9.2.2.

\m

{\bf Lemma 1.5.} {\it The family of distributions of\/ $\boldkey
H^\eps(\bullet)$ in the space\/ $\Bbb C[0,\, T]$ is tight
if the family of distributions of the initial points\/ $\boldkey
H^\eps(0)$ is tight\/} (in particular if the initial point
$\boldkey H^\eps(0)$ does not depend on $\,\eps$).

\s

{\bf Lemma 1.6.} {\it The family of distributions of\/ $\tilde{
\boldkey H}^\eps(\bullet)$ in the space\/ $\Bbb C[0,\, T]$ is tight
if the family of distributions of the initial points\/ $\tilde{
\boldkey H}^\eps(0)$ is tight\/} (in particular if the initial point
does not depend on $\,\eps$).

\m

The {\bf proofs} are by the use of Lemma 1.4 with $\,\beta = 4$,
$\alpha = 1$.

\m

The process $\tilde{\boldkey H}^\eps(t)$ is easier to handle than $\boldkey
H^\eps(t)$. One of the standard ways to prove that the distribution
of $\tilde{\boldkey H}^\eps(\bullet)$ converges weakly to a distribution
being the solution of a martingale problem is to establish that
$\tilde{\boldkey H}^\eps(t)$ is {\it approximately\/} (with the error of
the approximation going to $0$ as $\,\eps \to 0$) a solution
of that martingale problem; and that the solution of this problem
is unique. More specifically, we have to prove that for some
linear second-order differential operator $L$ and for some
(wide enough) class $\frak D$ of functions $\,f(\boldkey h)$
$$
    \EE\,\Bigl[f\bigl(\tilde{\boldkey H}^\eps(t)\bigr) - f(\boldkey h_0)
        - \int_0^t L f\bigl(\tilde{\boldkey H}^\eps(s)\bigr)\ ds\Bigr]
    \to 0
\eqno(1.19)
$$
as $\,\eps \to 0$, \ab uniformly for $\,t \in [0,\, T]$ (see [7],
Lemma 8.3.1).

\s

We can apply the It\^o formula to $\,f\bigl(\tilde{\boldkey H}^\eps(t)\bigr)$:
\ab for a twice continuously differentiable $\,f$, \ab using
(1.16), we get:
$$
\aligned
    &f\bigl(\tilde{\boldkey H}^\eps(t)\bigr) - f(\boldkey h_0)
    = \int_0^t \sum_{i = 1}^n \dfrac{\partial f}{\partial h_i}\ab
     \bigl(\tilde{\boldkey H}^\eps(s)\bigr) \cdot
     \nabla_{\boldkey w} u_i\bigl(\boldkey H^\eps(s),\,
     \boldsymbol\Phi^\eps(s),\,
        \boldkey W^\eps(s)\bigr)\ d\ab\boldkey W(s) \\
    &\ + \int_0^t \bigl[\sum_{i = 1}^n \dfrac{\partial f}{\partial h_i}\ab
     \bigl(\tilde{\boldkey H}^\eps(s)\bigr) \cdot
     \nabla_{\boldkey x}\ab u_i\bigl(\boldkey H^\eps(s),\,
     \boldsymbol\Phi^\eps(s),\, \boldkey W^\eps(s)\bigr)\cdot
     \boldkey b\ab\bigl(\boldkey H^\eps(s),\,
     \boldsymbol\Phi^\eps(s),\, \boldkey W^\eps(s)\bigr)            \\
    &\qquad+ \frac1{\ab2\ab} \sum_{i,\, j = 1}^n \dfrac{\partial^2 f}
        {\partial h_i\ab\partial h_j}\ab
        \bigl(\tilde{\boldkey H}^\eps(s)\bigr) \times \\
    &\qquad\qquad\ \ \times\nabla_{\boldkey w}\ab u_i\bigl(\boldkey H^\eps(s),\,
        \boldsymbol\Phi^\eps(s),\, \boldkey W^\eps(s)\bigr)\cdot
        \nabla_{\boldkey w}\ab u_j\bigl(\boldkey H^\eps(s),\,
            \boldsymbol\Phi^\eps(s),\, \boldkey W^\eps(s)\bigr)\bigr]\ ds\ab.
\endaligned
\eqno(1.20)
$$

If the random function under the sign of the stochastic integral
is square-integrable, the expectation of this integral is equal
to $0$, and we have something like formula (1.19), except
that we have an exact equality instead of ``$\to 0$'', and
that the factors by which the derivatives $\dfrac{\partial f}
{\partial h_i}\ab\bigl(\tilde{\boldkey H}^\eps(s)\bigr)$, $\dfrac{\partial^2 f}
{\partial h_i\ab\partial h_j}\ab\bigl(\tilde{\boldkey H}^\eps(s)\bigr)$
are multiplied are functions
not of $\tilde{\boldkey H}^\eps(s)$, but of $\boldkey H^\eps(s)$,
$\boldsymbol\Phi^\eps(s)$, $\boldkey W^\eps(s)$.

We can get rid of $\boldkey W^\eps(s)$ here by using (1.8) \ab:
we introduce the functions
$$
    \tilde B_i(\boldkey h,\, \boldsymbol\varphi)
    = \int_{\Bbb T^m} \nabla_{\boldkey x}\ab u_i(\boldkey h,
     \, \boldsymbol\varphi,\, \boldkey w) \cdot \boldkey b(\boldkey h,
     \, \boldsymbol\varphi,\, \boldkey w)\ d\boldkey w,
\eqno(1.21)
$$
$$
    \tilde A_{ij}(\boldkey h,\, \boldsymbol\varphi)
    = \int_{\Bbb T^m} \nabla_{\boldkey w}\ab u_i(\boldkey h,
     \, \boldsymbol\varphi,\, \boldkey w) \cdot \nabla_{\boldkey w}\ab u_j(
     \boldkey h,\, \boldsymbol\varphi,\, \boldkey w)\ d\boldkey w\ab
\eqno(1.22)
$$
(the integrands in (1.21) are twice differentiable in $(\boldkey
 h,\, \boldsymbol\varphi)$, and the integrands in (1.22)
three times), and get:
$$
\aligned
    \EE\ab\Bigl[f\bigl(\tilde{\boldkey H}^\eps(t)\bigr) - f(\boldkey h_0)\,
     - &\int_0^t \bigl[\dfrac1{\ab2\ab} \sum_{i,\, j = 1}^n
      \tilde A_{ij}\bigl(\boldkey H^\eps(s),\, \boldsymbol\Phi^\eps(s)\bigr)
      \cdot \dfrac{\partial^2 f}{\partial h_i\ab \partial h_j}\bigl(
      \tilde{\boldkey H}^\eps(s)\bigr)      \\
     &\ \ \ \ab + \sum_{i = 1}^n \tilde B_i\bigl(\boldkey H^\eps(s),\,
     \boldsymbol\Phi^\eps(s)\bigr)
      \cdot \dfrac{\partial f}{\partial h_i}\bigl(
      \tilde{\boldkey H}^\eps(s)\bigr)\bigr]\ ds\ab\Bigr] = O(\eps).
\endaligned
\eqno(1.23)
$$

Of course, $\boldkey H^\eps(s)$ is close to $\tilde{\boldkey H}^\eps(s)$,
so
$$
\aligned
    \EE\ab\Bigl[f\bigl(\boldkey H^\eps(t)\bigr) - f(\boldkey h_0)\,
     -  &\int_0^t \bigl[\dfrac1{\ab2\ab} \sum_{i,\, j = 1}^n
      \tilde A_{ij}\bigl(\boldkey H^\eps(s),\,
      \boldsymbol\Phi^\eps(s)\bigr)
      \cdot \dfrac{\partial^2 f}{\partial h_i\ab \partial h_j}\bigl(
      \boldkey H^\eps(s)\bigr)      \\
     &\ \ \ + \sum_{i = 1}^n \tilde B_i\bigl(\boldkey H^\eps(s),\,
     \boldsymbol\Phi^\eps(s)\bigr)
      \cdot \dfrac{\partial f}{\partial h_i}\bigl(
      \boldkey H^\eps(s)\bigr)\bigr]\ ds\ab\Bigr] = O(\eps)\ab.
\endaligned
\eqno(1.24)
$$

For the future use let's write the version of (1.24) with
$\,e^{- \lambda s}\ab$ and the integral from $0$ to $\infty$\ab:
$$
\aligned
    \EE  &
    \int_0^\infty e^{- \lambda t} \cdot\bigl[
    \dfrac1{\ab2\ab} \sum_{i,\, j = 1}^n
      \tilde A_{ij}\bigl(\boldkey H^\eps(s),\,
      \boldsymbol\Phi^\eps(s)\bigr)
      \cdot \dfrac{\partial^2 f}{\partial h_i\ab \partial h_j}\bigl(
      \boldkey H^\eps(s)\bigr)   \\
    &\ \ \ + \sum_{i = 1}^n \tilde B_i\bigl(\boldkey H^\eps(s),\,
     \boldsymbol\Phi^\eps(s)\bigr)
      \cdot \dfrac{\partial f}{\partial h_i}\bigl(
      \boldkey H^\eps(s)\bigr)
     - \lambda \cdot f\bigl(\boldkey H^\eps(t)\bigr)\bigr]\ ds
     = - f(\boldkey h_0) + O(\eps)\ab.
\endaligned
\eqno(1.25)
$$

Still we haven't reached (1,24): we have to get rid of
$\boldsymbol\Phi^\eps(s)$. Let us
average the functions $\ab\tilde A_{ik}(
\boldkey h,\, \boldsymbol\varphi)$, $\tilde B_i(\boldkey h,\,
\boldsymbol\varphi)$ over $\,\boldsymbol\varphi \in \Bbb T^p$:
$$
    \overline A_{i j}(\boldkey h)
    = \int_{\Bbb T^p} \tilde A_{ij}(\boldkey h,\,
    \boldsymbol\varphi)\ d\boldsymbol\varphi,
    \qquad \overline B_i(\boldkey h)
    = \int_{\Bbb T^p} \tilde B_i(\boldkey h,\,
    \boldsymbol\varphi)\ d\boldsymbol\varphi.
\eqno(1.26)
$$

We'll handle integrals of the type we have in (1.24), (1.25)
in the next section; but now we are finally in a position to formulate
our main result.

\s

Let us introduce two main conditions imposed on the system
(1.1). One of them has to do only with the unperturbed system
(0.7):

\s
\noi Condition $\ab\star$\,:

For every nonzero vector $\,\boldkey k = (k_1,\ab...,\, k_p)$ with
components being integers, let $N\!ull = N\!ull_{\boldkey k}
= \{\boldkey h\col \sum_{j = 1}^{\hm p} k_j \cdot \omega_j(
\boldkey h) = 0\}$. We require that the part of this set within
every bounded region should consist of finitely many points;
plus, outside of arbitrary small neighborhoods of these points,
finitely many smooth curves;
plus, outside of arbitrarily small neighborhoods of these
curves, finitely many smooth two-dimensional surfaces; ...;
and finally, outside arbitrarily small neighborhoods of the
mentioned $(n - 2)$-dimensional surfaces, a finite number
of smooth $(n - 1)$-dimensional surfaces (in the case of $\,
n = 1$, the part of this set within every finite interval
should just consist of finitely many points).

The set $\boldkey N\!\boldkey u\boldkey l\boldkey l
= \bigcup_{\boldkey k \neq \bold0} N\!ull_{\boldkey k}$ is the
set of $\,\boldkey h\,$ for which $\,\omega_1(\boldkey h)$,
..., $\omega_p(\boldkey h)\ab$ are rationally dependent (the
resonance set). We
mentioned in Section 0 the natural condition that the Lebesgue
measure $\,\lambda_n(\boldkey N\!\boldkey u\boldkey l\boldkey l)$
should be equal to $0$. Condition $\star$ is stronger
and more specific than the condition $\,\lambda_n(
\boldkey N\!\boldkey u\boldkey l\boldkey l) = 0$.\s

\s

The other condition, in contrast, has to do only with the
$\,\boldkey h$-component $\,\boldkey b_{\boldkey h}(\boldkey h,
\, \boldsymbol \varphi,\, \boldkey w)\ab$ of the perturbations:

\s

\noi Condition $\star\star$\ab:

The matrix\/ $\bigl(\tilde A_{ij}(\boldkey h,\, \boldsymbol
\varphi)\bigr)_{i,\ab j = 1,\ab...,\, n}$ defined by (1.22)
is uniformly positive definite\ab:
$$
    \sum_{i,\, j  = 1}^n \tilde A_{ij}(\boldkey h,\, \boldsymbol \varphi)
    \cdot \xi_i\, \xi_j
    \geq \underline a \cdot |\boldsymbol \xi|^2, \qquad \underline a > 0,
\eqno(1.27)
$$
for all\/ $\,\boldsymbol \xi \in \Bbb R^n,$ $\ab \boldkey h \in \Bbb R^n,$
$\boldsymbol\varphi \in \Bbb T^p.$

\m

{\bf Theorem 1.1.} {\it Let the function\/ $\,\boldkey b(\boldkey x,\,
\boldkey w) = \boldkey b(\boldkey h,\, \boldsymbol\varphi,\, \boldkey w)\ab$
be bounded, continuous in\/ $(\boldkey h,\, \boldsymbol\varphi,\, \boldkey w)$
and H\"older-continuous in\/ $\,\boldkey w,$ \ab uniformly in\/
$(\boldkey h,\, \boldsymbol\varphi),$ together with its derivatives
in\/ $\boldkey h,$ $\boldsymbol\varphi\,$ up to order\/~$3.$ \! Let\/
{\rm Conditions} $\star$ and\/ $\star\star$ be satisfied.

\s

Then\/ $\boldkey H^\eps(t)$ converges in distribution
to the diffusion process with the generating operator
$$
    L = \dfrac1{\ab2\ab} \sum_{i,\, j = 1}^n \overline A_{ij}(
     \boldkey h) \cdot \dfrac{\partial^2}{\partial h_i\ab\partial h_j}
    + \sum_{i = 1}^n \overline B_i(\boldkey h) \cdot
     \dfrac{\partial}{\partial h_i}
\eqno(1.28)
$$
and the same initial point\/}.

\m

In Section 2 we are going to prove Theorem 1.1.

\bigskip

\leftline{\bf 2. Averaging over the fast component. Proof of Theorem 1.1.}

\m

Let $N$ be the intersection of the set  $N\!ull = N\!ull_{\boldkey h}$
with a bounded set; for $\,\gamma > 0\ab$ let $N_{+\gamma}$
be its  $\,\gamma $-neighborhood

We want to estimate the expectation of the time spent by the
process $\boldkey H^\eps(t)$ in $N_{+\gamma}$ up to time~$T$.
This is approximately the same as estimating $\botsmash{\EE \dsize\int_0^\infty}
e^{- \lambda s}\, I_{N_{+ \gamma}}\bigl(\boldkey H^\eps(s)\bigr)
\ ds$ for some $\,\lambda > 0$: \ab indeed,
$$
    \EE \int_0^T I_{N_{+ \gamma}}\bigl(\boldkey H^\eps(s)\bigr)\ ds
    \leq e^{\lambda T} \cdot \EE \int_0^\infty e^{- \lambda s}\,
    I_{N_{+ \gamma}}\bigl(\boldkey H^\eps(s)\bigr)\ ds\ab.
\eqno(2.1)
$$

{\bf Lemma 2.1.} {\it Under the conditions of\/ {\rm Theorem 1.1}
there exists a constant\/ $C$ such that for the process starting
from any initial point}
$$
    \EE \int_0^\infty e^{- \lambda s}\, I_{N_{+ \gamma}}\bigl(
     \boldkey H^\eps(s)\bigr)\ ds
    \leq C \cdot \gamma + O(\eps).
\eqno(2.2)
$$

{\bf Proof.} Let us start with the one-dimensional case: $\ab n = 1$
(instead of the vector variable $\,\boldkey h\,$ and the vector-valued
stochastic process $\boldkey H^\eps(t)$ we'll have
the scalar coordinate\/~$\,h$, correspondingly, $H^\eps(t)$).
Let us start with the case of $N$ consisting of one point $\,h_*$.
Let $\,i(h)\ab$ be a continuous function that dominates the
indicator function $I_{S_{+\gamma}}(h)$: $\ab i(h) = 1\ab$
for $\,h \in [h_* - \gamma,\, h_* + \gamma]$, between $0$ and~
$1$ for $\,h \in [h_* - 2\ab\gamma,\, h_* - \gamma] \cup
[h_* + \gamma,\, h_* + 2\ab \gamma]$, and $\,i(h) = 0\ab$ outside
$[h_* - 2\ab \gamma,\, h_* + 2\ab \gamma]$. For some $\,\mu
> 0\ab$ let us consider the bounded solution $\,v(h)\ab$ of
the equation $\,\mu \cdot v(h) - \tfrac12\, v''(h) = i(h)$; it
is given by
$$
    v(h) = \int_{- \infty}^h \dfrac1{\sqrt{2\ab\mu}\,}\,
    e^{\sqrt{2\mu}\ab(\eta - h)} \cdot i(\eta)\ d\eta
    + \int_h^\infty \dfrac1{\sqrt{2\ab\mu}\,}\,
    e^{\sqrt{2\mu}\ab(h - \eta)} \cdot i(\eta)\ d\eta\ab.
\eqno(2.3)
$$
We have:
$$
    v'(h) = - \int_{- \infty}^h e^{\sqrt{2\mu}\ab(\eta - h)}
    \cdot i(\eta)\ d\eta
    + \int_h^\infty e^{\sqrt{2\mu}\ab(h - \eta)},
    \cdot i(\eta)\ d\eta\ab;
\eqno(2.4)
$$
$$
    v(h_0) \leq 4\ab\gamma/\sqrt{2\ab\mu},
\eqno(2.5)
$$
where $\,h_0\ab$ is the initial point $H^\eps(0)$.

From (2.4) we see that $|v'(h)| \leq \sqrt{2\ab\mu} \cdot v(h)$.

\s

We are going to write some inequality for the function in
brackets in formula (1.25) for the function $\,f = v$.
\ab For all $\,h$, $\boldsymbol \varphi\,$ we have:
$$
\aligned
    &\tfrac12\,\tilde A(h,\, \boldsymbol\varphi) \cdot
    v''(h) + \tilde B(h,\, \boldsymbol\varphi) \cdot v'(h)
    -\lambda \cdot v(h) \\
    &\quad= \,\tilde A(h,\, \boldsymbol\varphi) \cdot \mu \cdot v(h)
    - \tilde A(h,\, \boldsymbol\varphi) \cdot i(h)
    + \tilde B(h,\, \boldsymbol\varphi) \cdot v'(h)
    -\lambda \cdot v(h)     \\
    &\qquad\qquad\leq \overline a \cdot \mu \cdot v(h)
    - \underline a \cdot i(h)
    + \overline b \cdot \sqrt{2\ab\mu} \cdot v(h)
    - \lambda \cdot v(h).
\endaligned
\eqno(2.6)
$$
If we choose $\,\mu > 0\ab$ so small that $\,\overline  a
\cdot \mu \leq \lambda/2$, $\overline b \cdot \sqrt{2\ab\mu} \leq
\lambda/2$, \ab where $\,\overline a$, $\overline b\,$ are
the constants dominating $\tilde A_{11}(h,\, \boldsymbol
\varphi)$, $\tilde B_1(h,\, \boldsymbol\varphi)$, the expression
(2.5) is $\leq-\, i(h)$.

From this and formula (1.25) we get:
$$
    \EE \int_0^\infty e^{- \lambda s}\, i\bigl(G^\eps(s)\bigr)\ ds
    \leq \underline a^{- 1} \cdot v(h_0) + O(\eps);
\eqno(2.7)
$$
that is, in the one-dimensional case we have (2.2) satisfied
with $C = 4\, \underline a^{- 1}/\sqrt{2\ab\mu}\ab$.

Of course if the set $N$ consists of finitely many points, we
just multiply the right-hand side of (2.7) by the number of
these points.

\s

Now we go to the multidimensional case. Suppose $N$ is an
$(n - 1)$-dimensional surface described by the equation
$h_i = \!f(h_1,\ab...,\ab h_{i - 1}, h_{i + 1}, ...,\, h_n)$,
\ab where $\ab(h_1,...,\ab h_{i - 1}, h_{i + 1},...,
h_n)$ changes in a bounded region, and the function
$\,f\,$ is twice differentiable with bounded and continuous
derivatives. Let us extend the function $\,f(h_1,\ab...,\,
h_{i - 1},$ $h_{i + 1},\ab...,\, h_n)\ab$ from the bounded
region to the whole $\Bbb R^{n - 1}$ so that its first and
second derivatives are bounded an continuous. Let
us introduce in $\Bbb R^n$ a new coordinate system $(\tilde h_1,
\ab...,\, \tilde h_n)$ with\linebreak $\,\tilde h_i = h_i - f(h_1,\ab...,
\, h_{i - 1},\, h_{i + 1},\ab...,\, h_n)\ab$ and $\,\tilde h_j
= h_j\ab$ for $\,j \neq i$. \ab The differential operator
$\tfrac12 \sum_{j,\, l = 1}^{\te n} \tilde A_{jl}(\boldkey h,
\, \boldsymbol \varphi) \cdot \dfrac{\partial^2}{\partial h_j\,
\partial h_l} + \sum_{j = 1}^{\,n} \tilde B_j(\boldkey h,\,
\boldsymbol\varphi) \cdot \dfrac{\partial}{\partial h_j}$ is
written in the new coordinates $\,\tilde h_1$, ..., $\tilde h_n\ab$
as $\topsmash{\tfrac12 \sum_{j,\, l = 1}^{\te n}
\tilde{\tilde A\,}_{\!jl}(\boldkey h,
\, \boldsymbol \varphi) \cdot \dfrac{\partial^2}{\partial \tilde h_j\,
\partial \tilde h_l} + \sum_{j = 1}^{\,n} \tilde{\tilde B\ab }_j(\boldkey h,\,
\boldsymbol\varphi) \cdot \dfrac{\partial}{\partial \tilde h_j}}$;
the coefficients are bounded: $|\tilde{\tilde A\,}_{\!jl}(\boldkey h,
\, \boldsymbol \varphi)| \leq \tilde{\overline a}$, $|\tilde{\tilde B\ab}_j
(\boldkey h,\, \boldsymbol \varphi)| \leq \tilde{\overline b}$,
and
$$
    \tilde{\tilde A\,}_{\!ii}(\boldkey h,\, \boldsymbol \varphi)
    = \sum_{j,\, l = 1}^n \xi_j\,\xi_l
    \geq \underline a \cdot \sum_{j = 1}^n \xi_j^2\ab,
\eqno(2.8)
$$
where $\,\xi_i = 1$, $\xi_j = \topsmash{\dfrac{\partial f}{\partial h_j}\ab}$
for $\,j \neq i$. \ab So we have $\tilde{\tilde A\,}_{\!ii}(
\boldkey h,\, \boldsymbol \varphi) \geq \underline a$ ($> 0)$.

The $\,\gamma$-neighborhood \ab $N_{+ \gamma}$ lies within
the set $\{\boldkey h\col |\tilde h_i| \leq c \cdot \gamma\}$,
where $\,c =$\linebreak $\sup \sqrt{1 + \sum_{j \neq i}
(\partial f/\partial h_j)^2}$. Take a function
$\,i(\tilde h_i)\ab$ dominating the indicator function \linebreak
$I_{[- c\,\gamma,\, c\,\gamma]}(\tilde h_i)$ and equal to $0$
outside the interval $[- 2\ab c\ab \gamma,\, 2\ab c\ab \gamma]$;
and consider the function $\,V(\boldkey h) = v(\tilde h_i)$
given by (2.3) with $\,h = h_i$. We have the estimate
(2.7), and so (2.2) with $C  = 4\ab c\ab
\underline a^{- 1}/\sqrt{2\ab \mu}\ab$. Because of (2.1),
we have the same kind of estimate for $\EE \dsize\int_0^T I_{N_{+ \gamma}}
\bigl(\boldkey H^\eps(s)\bigr)\ ds$\ab: we just multiply it by
$\,e^{\lambda T}$.

We take care of surfaces with finitely many pieces described
by different equations by adding the corresponding estimates;
abd of surfaces of smaller dimensions by noticing that a
$\,\gamma$-neighborhood \ab of an $\,s$-dimensional, $\ab s < n -1$,
smooth surface is a part of the union of $\,\gamma$-neighborhoods
of finitely many $(n - 1)$-dimensional smooth surfaces.

\bigskip

{\bf Lemma 2.2.} {\it Let the conditions of\/ {\rm Theorem 1.1}
be satisfied$\ab.$ Let\/ $\,g(\boldkey h,\, \boldsymbol\varphi)\ab$ be
a bounded uniformly continuous function$\ab;$ $\ab\overline g(\boldkey h)
= \dsize\int_{\Bbb T^p} g(\boldkey h,\, \boldsymbol\varphi)\
d\boldsymbol\varphi\ab.$

Then for every\/ $T > 0$ we have\ab}:
$$
    \EE\ab\max_{0 \leq t \leq T}\ab
    \Bigl|\int_0^t g\bigl({\boldkey H}^\eps(s),\,
    \boldsymbol\Phi^\eps(s)\bigr)\ ds
    - \int_0^t \overline g\bigl({\boldkey H}^\eps(s)\bigr)\ ds\ab\Bigr|
    \to 0 \quad\ (\eps \to 0).
\eqno(2.9)
$$

{\bf Proof.} We need to prove that for every $\,\delta > 0\ab$
for sufficiently small $\,\eps\,$
$$
    \EE\ab\max_{0 \leq t \leq T}\ab
    \Bigl|\int_0^t g\bigl({\boldkey H}^\eps(s),\,
    \boldsymbol\Phi^\eps(s)\bigr)\ ds
    - \int_0^t \overline g\bigl({\boldkey H}^\eps(s)\bigr)\ ds\ab\Bigr|
    < \delta.
\eqno(2.10)
$$
Let us approximate the function $\,g(\boldkey h,\,
\boldkey z) - \overline g(\boldkey h)\ab$ up to $\,\delta/2\ab T\ab$
in the $\Bbb C(\Bbb R^n \times \Bbb T^p)$-norm by a finite
trigonometric sum
$$
    \sum_{\boldkey k\ab =\ab (k_1,\,...,\, k_p) \neq \,\bold 0, \
    - K \leq k_j \leq K}
    G_{\boldkey k}(\boldkey h) \cdot \exp\bigl\{2\ab\pi\ab i \cdot
     \sum_{j = 1}^p k_j \cdot \varphi_j\bigr\},
\eqno(2.11)
$$
where the functions $C_{\boldkey k}(\boldkey h)$ are bounded
and continuous with their second derivatives. It's enough
to consider the functions
$$
    g_{\boldkey k}(\boldkey h,\, \boldsymbol\varphi)
    = C_{\boldkey k}(\boldkey h) \cdot \exp\bigl\{2\ab\pi\ab i \cdot
     \sum_{j = 1}^p k_j \cdot \varphi_j\bigr\}, \qquad
     \boldkey k \neq \bold 0
\eqno(2.12)
$$
(note that $\,\overline g_{\boldkey k}(
\boldkey h) = 0$). Let us introduce the random variables
$$
    M\!ax^{\,\eps}_{\boldkey k}
    = \max_{0 \leq t \leq T}\ab\Bigl|\int_0^t
     g_{\boldkey k}\bigl(\boldkey H^\eps(s),\,
     \boldsymbol\Phi^\eps(s)\bigr)\ ds\ab\Bigr|\ab.
\eqno(2.13)
$$
It's enough to prove that for sufficiently small $\,\eps$
$$
    \EE\, M\!ax^{\,\eps}_{\boldkey k} < \delta'
    = \dfrac\delta{2\ab T\ab (2\ab K)^p}.
\eqno(2.14)
$$

Let us consider the event
$$
    A_R = \bigl\{\max_{0 \leq t \leq T}\ab |\boldkey H^\eps(t)| \leq R\bigr\}
\eqno(2.15)
$$
and its complement $A_R^{\,c}$. For sufficiently large $R$ we have
$$
    \PP\bigl\{\max_{0 \leq t \leq T}\ab |\tilde{\boldkey H}^\eps(t)|
    > R - 1\bigr\} \leq \dfrac{\delta'}{3\ab\|g_{\boldkey k}\|},
\eqno(2.16)
$$
(the random function $\tilde{\boldkey H}^\eps(t)$ is described by
formula (1.16) with an integral and a stochastic integral, and
we obtain (2.16) using the Kolmogorov inequality); and since
the distance between $\boldkey H^\eps(t)$ and $\tilde{\boldkey H}^\eps(t)$
is $O(\eps)$, for sufficiently small $\,\eps\,$ we have:
$$
    \botsmash{\PP(A_R^{\,c})
    \leq \PP\bigl\{\max_{0 \leq t \leq T}\ab |\tilde{\boldkey H}^\eps(t)|
    > R - 1\bigr\} \leq \dfrac{\delta'}{3\ab\|g_{\boldkey k}\|},}
\eqno(2.17)
$$
and
$$
    \EE\, M\!ax^{\,\eps}_{\boldkey k}
    = \EE\ab [I_{A_R} \cdot M\!ax^{\,\eps}_{\boldkey k}]
     + \EE\ab [I_{A_R^{\,c}} \cdot M\!ax^{\,\eps}_{\boldkey k}]
     \leq \EE\ab [I_{A_R} \cdot M\!ax^{\,\eps}_{\boldkey k}]
     + \delta'\!/3.
\eqno(2.18)
$$

Now let us choose a $\Delta = \Delta(\eps)$ such that $M =
T/\Delta$ is an integer, $\Delta \to 0$, $\Delta/\eps \to \infty$
($\eps \to 0$). Let $\,t_0 \in [0, \Delta)$, $t_i = t_0 + i
\cdot \Delta$ (the choice of $\,t_0\ab$ will
be specified later). The expectation of $I_{A_R} \cdot
M\!ax^{\,\eps}_{\boldkey k}$ does not exceed
$$
\aligned
    &\EE\ab\Bigl[I_{A_R} \cdot \max_{0 \leq t \leq t_0}\ab
     \Bigl|\int_0^t g_{\boldkey k}\bigl(\boldkey H^\eps(s),\,
     \boldsymbol\Phi^\eps(s)\bigr)\ ds\ab\Bigr|\ab\Bigr] \\
     &\qquad\ + \sum_{i = 0}^{M - 2} \EE\ab\Bigl[I_{A_R} \cdot
     \max_{t_i \leq t \leq t_{i + 1}}\ab
     \Bigl|\int_{t_i}^t g_{\boldkey k}\bigl(\boldkey H^\eps(s),\,
     \boldsymbol\Phi^\eps(s)\bigr)\ ds\ab\Bigr|\ab\Bigr] \\
     &\qquad\qquad\qquad+ \EE\ab\Bigl[I_{A_R} \cdot
     \max_{T - \Delta + t_0 \leq t \leq T}\ab
     \Bigl|\int_{t_i}^t g_{\boldkey k}\bigl(\boldkey H^\eps(s),\,
     \boldsymbol\Phi^\eps(s)\bigr)\ ds\ab\Bigr|\ab\Bigr]\ab.
\endaligned
\eqno(2.19)
$$
The first and the last summands are $O(\Delta) \to 0$; of course
the same is true for all other summands, but the number of
summands is approximately $T/\Delta$, which gives us only that
the total expectation is $O(1)$\ab: not enough. We can do better.

The stochastic process $\bigl(\boldkey H^\eps(t),\, \boldsymbol
\Phi^\eps(t),\, \boldkey W^\eps(t)\bigr)$ is a time-homogeneous Markov one;
we'll be denoting $\EE_{\boldkey h,\, \boldsymbol\varphi,\, \boldkey w}$
the expectation associated with this process evaluated under
the assumption that $\boldkey H^\eps(0) = \boldkey h$, $\boldsymbol\Phi^\eps(0)
= \boldsymbol\varphi$, $\boldkey W^\eps(0) = \boldkey w$ (this
means, in particular, that we take $(\boldkey h,\, \boldsymbol\varphi)$
as the initial point $\,\boldkey x_0$; \ab the expectation
that we denoted with just $\EE$ being associated with an arbitrary
initial distribution).  The Markov
property with respect to the time $\,t_i\ab$ yields for the
$\,i$-th summand in (2.19) (with the exception of the last
summand, which is $O(\Delta)$ anyway):
$$
    \EE\ab\bigl[I_{\{\max_{\,0 \leq t \leq t_i} |\boldkey H^\eps(t)| \leq R\}}
    \cdot M^\eps_i\bigl(\boldkey H^\eps(t_i),\, \boldsymbol\Phi^\eps(t_i),\,
    \boldkey W^\eps(t_i)\bigr)\bigr],
\eqno(2.20)
$$
where
$$
    M^\eps_i(\boldkey h,\, \boldsymbol\varphi,\, \boldkey w)
    = \EE_{\boldkey h,\, \boldsymbol\varphi,\, \boldkey w} \Bigl[
    I_{\{\max_{\,0 \leq t \leq T - t_i} |\boldkey H^\eps(t)| \leq R\}}
    \cdot \max_{0 \leq t \leq \Delta} \Big|\int_0^t g_{\boldkey k}\bigl(
     \boldkey H^\eps(s),\, \boldsymbol\Phi^\eps(s)\bigr)\ ds\ab\Bigr|\ab\Bigr].
\eqno(2.21)
$$

We have $\,\max_{\,0 \leq t \leq T- t_i} |\boldkey H^\eps(t)
- \tilde{\boldkey H}^\eps(t)| = O(\eps)$, $\max_{\,0 \leq t \leq T- t_i}
|\boldsymbol\Phi^\eps(t) - \tilde{\boldsymbol\Phi}^\eps(t)|$
$= O(\eps)$, \ab and for all $\, t \in [0,\, \Delta]$
$$
    \int_0^t g_{\boldkey k}\bigl(\boldkey H^\eps(s),\,
     \boldsymbol\Phi^\eps(s)\bigr)\ ds
    - \int_0^t g_{\boldkey k}\bigl(\tilde{\boldkey H}^\eps(s),\,
     \tilde{\boldsymbol\Phi}^\eps(s)\bigr)\ ds
    = O(t \cdot \eps) = O(\Delta \cdot \eps);
\eqno(2.22)
$$
so the expectation (2.21) is not greater than
$$
    \EE_{\boldkey h,\, \boldsymbol\varphi,\, \boldkey w} \Bigl[
    I_{\{\max_{\,0 \leq t \leq T - t_i} |\tilde{\boldkey H}^\eps(t)| \leq R + 1\}}
    \cdot \max_{0 \leq t \leq \Delta} \Bigl|\int_0^t g_{\boldkey k}\bigl(
     \tilde{\boldkey H}^\eps(s),\,
     \tilde{\boldsymbol\Phi}^\eps(s)\bigr)\ ds\ab\Bigr|\ab\Bigr]
    + O(\Delta \cdot \eps).
\eqno(2.23)
$$

Because of the factor $I_{\{\max_{\,0 \leq t \leq t_i} |\boldkey H^\eps(t)|
\leq R\}}$ in formula (2.20) we need to estimate the expectation
in (2.23) only for $(\boldkey h,\, \boldsymbol\varphi,\, \boldkey w)$
with~$|\boldkey h| \leq R$.

Note that we apply to processes $\bigl(\boldkey H^\eps(t),\,
\boldsymbol\Phi^\eps(t),\, \boldkey W^\eps(t)\bigr)$ and $\bigl(\tilde{
\boldkey H}^\eps(t),\, \tilde{\boldsymbol\Phi}^\eps(t)\bigr)$ different
methods: Markov-process methods; and stochastic-integrals methods,
correspondingly. We are going to use essentially the same technique as
we used in Lemma 1.2: in both cases we introduce an auxiliary
function ($U(\boldkey x,\, \boldkey w)$ in one case, $v(\boldkey h,\,
\boldkey z)$ in the other) being a solution of a partial differential
equation and apply the It\^o formula. Our present case is simpler
only in that the solution is written as a simple explicit formula,
but it is more complicated due to the fact that the equation
is not satisfied on the whole space. We have to circumvent it.

\s

Let $\,v(\boldkey h,\, \boldsymbol\varphi)\ab$ be a function on $\Bbb R^n
\times \Bbb T^p$ that is bounded and continuous together with its
second derivatives. Let us apply the It\^o formula to $\,\eps \cdot
v\bigl(\tilde{\boldkey H}^\eps(t),\, \tilde{\boldsymbol\Phi}^\eps(t)\bigr)\ab$:
$$
\aligned
    &\int_0^t \nabla_{\boldsymbol\varphi}\, v\bigl(
     \tilde{\boldkey H}^\eps(s),\, \tilde{\boldsymbol\Phi}^\eps(s)\bigr)
     \cdot \boldsymbol\omega\bigl(\boldkey H^\eps(s)\bigr)\ ds \\
    &\ \ \ = \eps \cdot \bigl[v\bigl(\tilde{\boldkey H}^\eps(t),\,
     \tilde{\boldsymbol\Phi}^\eps(t)\bigr)
    - v\bigl(\boldkey H^\eps(0),\, \boldsymbol\Phi^\eps(0)\bigr)\bigr] \\
    &\qquad\ \ \ + \text{\,some integrals and stochastic integrals with integrands
     being } O(\eps)\ab,
\endaligned
\eqno(2.24)
$$
and also we can replace $\,\boldsymbol\omega\bigl(\boldkey H^\eps(s)\bigr)\ab$
with $\,\boldsymbol\omega\bigl(\tilde{\boldkey H}^\eps(s)\bigr)$:
$$
\aligned
    &\int_0^t \boldsymbol\omega\bigl(\tilde{\boldkey H}^\eps(s)\bigr)
    \cdot \nabla_{\boldsymbol\varphi}\, v\bigl(\tilde{\boldkey H}^\eps(s),\,
     \tilde{\boldsymbol\Phi}^\eps(s)\bigr)\ ds \\
    &\ \ \ \ = O(\eps) + \text{\,some integrals and stochastic
    integrals with integrands being } O(\eps)
\endaligned
\eqno(2.25)
$$
(another integral with integrand of order $O(\eps)$ will be added).

\s

If we were able to produce a function $\,v(\boldkey h,\,
\boldsymbol\varphi)\ab$ being a solution of the equation
$$
    \boldsymbol\omega(\boldkey h) \cdot \nabla_{\boldsymbol\varphi}\ab
    v(\boldkey h,\, \boldsymbol\varphi)
    = g_{\boldkey k}(\boldkey h,\, \boldsymbol\varphi),
\eqno(2.26)
$$
Lemma 2.2 would have been proved. The function
$$
    v(\boldkey h,\, \boldsymbol\varphi)
    = \dfrac{G_{\boldkey k}(\boldkey h) \cdot
    \exp\bigl\{2\ab\pi\ab i \cdot
    \sum_{j = 1}^{\ p} k_j \cdot \varphi_j\bigr\}}
    {2\ab\pi\ab i \cdot \sum_{j = 1}^{\ p} k_j \cdot \omega_j(\boldkey h)};
\eqno(2.27)
$$
(not defined everywhere!) is such a solution on the set of its
arguments for which the denominator is not equal to $0$. Let
us modify this function so it is defined on the whole space
and smooth on it.

\s

Let us choose a positive $\,\gamma$ (we are going to specify
its choice later). Let $\, v(\boldkey h,\, \boldsymbol\varphi)\ab$
be given by formula (2.27) for $\,\boldkey h \in N\!ull_{+\,\gamma}^{\ c}$,
$|\boldkey h| \leq R + 1$ (bounded and continuous together
with its second derivatives on the set $\bigl(N\!ull_{+\,\gamma}^{\ c}
\cap \{\boldkey h\col |\boldkey h| \leq R + 1\}\bigr)
\times \Bbb T ^p$), and let us extend it to the whole $\Bbb R^n
\times \Bbb T^p$ as a function that is bounded and continuous
together with its second derivatives..

\s

Since $\,\boldsymbol\omega(\boldkey h) \cdot \nabla_{\boldsymbol\varphi}\,
v(\boldkey h,\, \boldsymbol\varphi) = g_{\boldkey k}(\boldkey h,\,
\boldsymbol\varphi)$ \ab
for all $\,\boldkey h \in N\!ull_{+\,\gamma}^{\ c}$\ab, $|\boldkey h|
\leq R$, by (2.19) we have
$$
\aligned
    &\int_0^t g_{\boldkey k}\bigl(\tilde{\boldkey H}^\eps(s),\,
    \tilde{\boldsymbol\Phi}^\eps(s)\bigr)\ ds  \\
    &\ \ \ \ \, = O(\eps) + \text{\,some integrals and stochastic integrals
    with integrands being } O(\eps)
\endaligned
\eqno(2.28)
$$
as long as $|\tilde{\boldkey H}^\eps(s)| \leq R$, $\tilde{\boldkey H}^\eps(s)
\in N\!ull_{+\,\gamma}^{\ c}$ for all $\,s \in [0,\, t]$\ab.

\s

We'll use two different estimates for the expectation (2.23):
for $\,\boldkey h \in N\!ull_{+\,2\gamma}\ab$ it will be just
$O(\Delta)$, and the expectation (2.21) will also be $O(\Delta)$.
For $\,\boldkey h \in N\!ull^{\ c}_{+\,2\gamma}\ab$ the expectation
(2.23) is not greater than
$$
\aligned
    &\EE_{\boldkey h,\, \boldsymbol\varphi,\, \boldkey w} \Bigl[
    I_{\{\tilde{\boldkey H}^\eps(t) \in N\!ull_{+\,\gamma} \text{ or }
    |\tilde{\boldkey H}^\eps(t)| > R + 1 \text{ for some }
         \,t \in [0,\, \Delta]\}} \cdot
    \max_{0 \leq t \leq \Delta} \Bigl|\int_0^t g_{\boldkey k}\bigl(
     \tilde{\boldkey H}^\eps(s),\,
     \tilde{\boldsymbol\Phi}^\eps(s)\bigr)\ ds\ab\Bigr|\ab\Bigr]\\
    &\ + \EE_{\boldkey h,\, \boldsymbol\varphi,\, \boldkey w} \Bigl[
    I_{\{\tilde{\boldkey H}^\eps(t) \in N\!ull^{\ c}_{+\,\gamma}\,,\,
    |\tilde{\boldkey H}^\eps(t)| \leq R + 1 \text{ for all }
         \,t \in [0,\, \Delta]\}} \cdot
    \max_{0 \leq t \leq \Delta} \Bigl|\int_0^t g_{\boldkey k}\bigl(
     \tilde{\boldkey H}^\eps(s),\,
     \tilde{\boldsymbol\Phi}^\eps(s)\bigr)\ ds\ab\Bigr|\ab\Bigr].
\endaligned
\eqno(2.29)
$$
The first summand is not greater than
$$
    \PP_{\boldkey h,\, \boldsymbol\varphi,\, \boldkey w}\ab
     \{\tilde{\boldkey H}^\eps(t) \in N\!ull_{+\,\gamma} \text{ or }
    |\tilde{\boldkey H}^\eps(t)| > R + 1 \text{ for all }
         \,t \in [0,\, \Delta]\} \cdot O(\Delta)\ab;
\eqno(2.30)
$$
For $(\boldkey h,\, \boldsymbol\varphi,\, \boldkey w)$ with $\,\boldkey h
\in N\!ull^{\ c}_{+\,2\gamma}$\ab, $|\boldkey h| \leq R$ the
probability in (2.30) goes to $0$, uniformly for all mentioned
points $(\boldkey h,\, \boldsymbol\varphi,\, \boldkey w)$ (proved using
the Kolmogorov inequality: remember that the time interval goes
to $0$), and the first summand in (2.29) is $\,o(\Delta)$.

In the second expectation in (2.29) we can use formula (2.20),
and it is $O(\eps)$; since $\,\eps = o(\Delta)$, \ab the expectation
(2.23), and (2.21) with it, is $\,o(\Delta)\ab$ for $\,\boldkey h
\in N\!ull_{+ 2\gamma}^{\ c}$, $|\boldkey h| \leq R$.

So the sum (2.19) is not greater than
$$
\aligned
    &\sum_{i = 0}^{M - 2}
    \PP\{\boldkey H^\eps(t_i)
    \in N\!ull_{+ 2 \gamma}^{\ c}\ab,\, |\boldkey H^\eps(t_i)| \leq R\}
     \cdot o(\Delta)        \\
    &\qquad\qquad\quad+ \sum_{i = 0}^{M - 2} \PP\{\boldkey H^\eps(t_i)
    \in N\!ull_{+ 2 \gamma}\ab,\, |\boldkey H^\eps(t_i)| \leq R\}
     \cdot O(\Delta) + O(\Delta).
\endaligned
\eqno(2.31)
$$
The first term here is not greater than the number of summands
(which is $T/\Delta - 1$), multiplied by $o(\Delta)$; so
the first term is $o(1)$. The second term in (2.31) is not
greater than
$$
    \sum_{i = 0}^{M - 2} \PP\{\tilde{\boldkey H}^\eps(t_i)
    \in N\!ull_{+ 3 \gamma}\ab,\,
    |\tilde{\boldkey H}^\eps(t_i)| \leq R + 1\}
     \cdot O(\Delta).
\eqno(2.32)
$$
The sum here is equal to
$
    \,\EE \sum_{\,i = 0}^{M - 2} I_{\{\tilde{\boldkey H}^\eps(t_0 +\, i \cdot \Delta)
    \in N\!ull_{+ 3 \gamma}\ab,\, |\tilde{\boldkey H}^\eps(t_i)| \leq R + 1\}}.
$
The average of this sum over $\,t_0 \in [0,\, \Delta)\ab$ is
$$
\aligned
    &\dfrac1{\ab\Delta\ab} \cdot \EE \int_0^\Delta
    \sum_{i = 0}^{M - 2} I_{\{\tilde{\boldkey H}^\eps(t_0 +\, i \cdot \Delta) \in
    N\!ull_{+ 3 \gamma}\ab,\, |\tilde{\boldkey H}^\eps(t_0 + i \cdot \Delta)|
    \leq R + 1\}}\ dt_0 \\
    &\qquad\quad = \dfrac1{\ab\Delta\ab} \cdot \EE \sum_{i = 0}^{M - 2}
    \int_{i \cdot \Delta}^{(i + 1) \cdot \Delta}
    I_{\{\tilde{\boldkey H}^\eps(t) \in
    N\!ull_{+ 3 \gamma}\ab,\,
    |\tilde{\boldkey H}^\eps(t)| \leq R + 1\}}\ dt        \\
    &\qquad\qquad\qquad\qquad\qquad\quad
    = \dfrac1{\ab\Delta\ab} \cdot \EE \int_0^{T - \Delta}
    I_{\{\tilde{\boldkey H}^\eps(t) \in
    N\!ull_{+ 3 \gamma}\ab,\,
    |\tilde{\boldkey H}^\eps(t)| \leq R + 1\}}\ dt.
\endaligned
\eqno(2.33)
$$
There exists at least one $\,t_0 \in [0,\, \Delta)\ab$ for
which the sum in (2.32) is not greater than its average;
this $\,t_0\ab$ will be our choice (see above).

So by Lemma 2.1 with $\ab3\ab \gamma\,$ instead of $\,\gamma\,$
we have for $\,t_i = t_0 + i\ab \Delta\ab$ with the above
choice of $\,t_0$:
$$
    \dfrac1\Delta \sum_{i = 0}^{M - 2} \PP\{
    \tilde{\boldkey H}^\eps(t_i) \in N\!ull_{+ 3\ab\gamma},\,
    |\tilde{\boldkey H}^\eps(t_i)| \leq R + 1\}
    \leq \dfrac1\Delta \cdot [3\ab C\ab \gamma + O(\eps)];
\eqno(2.34)
$$
by choosing $\,\gamma\,$ and $\,\eps\,$ small enough we achieve
the inequality (2.14), which proves\linebreak Lemma~2.2.

\bigskip

Let us return to formula (1.19). As the set $\frak D$ we take
the set of all functions that are bounded and continuous together
with their second derivatives; and (1.19)
follows from formula (1.24) and Lemma 2.1 (the functions
$\tilde A_{ij}(\boldkey h,\, \boldsymbol \varphi) \cdot\! \botsmash{
\dfrac{\partial^2 f}{\partial h_i\ab \partial h_j}(\boldkey h)}$,
$\tilde B_i(\boldkey h,\, \boldsymbol \varphi) \cdot \botsmash{\dfrac{
\partial f}{\partial h_i}}(\boldkey h)$ \ab are bounded and continuous).

\s

{\bf Lemma 2.3.} {\it Let formula\/ $(1.19)$ hold for\/ $\,f \in
\frak D;$ let the solution of the martingale problem associated
with the operator\/ $L$ considered on\/ $\frak D,$ with a prescribed
initial distribution$,$ be unique$.$ The the distribution of\/
$\tilde{\boldkey H}^ \eps(\bullet)$ converges weakly to the solution
of this martingale problem\ab}.

\s

{\bf Proof.} Remember that the family of distributions of $\tilde{
\boldkey H}^\eps(\bullet)$ is tight (Lemma 1.5); then it is
the standard reasoning for proving weak convergence.

\m

{\bf Lemma 2.4.} {\it Let for every\/ $\,\lambda > 0\ab$ and
for every\/ $\,f\,$ belonging to a distinguishing set\/~$\frak H$
{\rm(that is, such a set that for finite measures $\,\mu_1$,
$\mu_2\ab$ the coincidence $\dsize\int f\ d\mu_1 = \dsize\int
f\ d\mu_2$ for all $\,f \in \frak H\ab\ab$ implies $\,\mu_1 = \mu_2$;
forgot the usual term for it)} there is a solution\/ $F \in
\frak D$ of the equation\/ $L\ab F - \lambda\ab F = f.$ Then
the solution of the martingale problem associated with\/ $L$
with a prescribed initial distribution is unique\ab}.

\s

See [7], Lemma 8.3.1.

\m

The following Lemma is a standard fact from the theory of partial
differential equations:

\s

{\bf Lemma 2.5.} {\it Let the functions\/ $\,a_{ij}(\boldkey h),$
$b_i(\boldkey h),$ $\boldkey h \in \Bbb R^n,$ \ab be bounded
and continuous together with their second derivatives\/$,$ and
let the elliptic operator\/ 
$L = \botsmash{\tfrac1{\ab2\ab} \sum_{\,i,\, j} a_{ij}(\boldkey h)} \times $
$\times \dfrac{\partial^2}{\partial h_i\, \partial h_j} + \sum_{\,i}
b_i(\boldkey h) \cdot \dfrac{\partial}{\partial h_i}$ be 
uniformly non-degenerate$.$
Then a solution \/ $F \in \frak D$ of the equation\/ $\ab L\ab F
- \lambda\ab F = f\,$ exists for all\/ $\,\lambda > 0\ab$ and all
right-hand sides\/ $\,f\,$ belonging to the set\/~$\frak H$
of infinitely differentiable functions that are equal to\/~$0$
outside a compact set\/} (different for different $\,f \in
\frak H$). 

\s

This proves Theorem 1.1.

\s

Condition $\star$ of this theorem is rather complicated, involving
surfaces of different dimensions and separate points; this
may be necessary
if the functions $\sum k_j \cdot \omega_j(\boldkey h)$ have
critical points.

In the case $\,n = 1\ab$ the condition of Theorem 1.1 about
finitely many points in $N\!ull_{\boldkey k}$ in every bounded
region can be replaced by the condition that the Lebesgue measure
of the set of $\,h\,$ for which $\,\omega_1(h)$, ..., $\omega_p(h)\ab$
are rationally dependent is equal to $0$. This is because for
sufficiently small $\,\gamma\,$ the closed set $N\!ull_{\boldkey k}
\cap [- R,\, R]$ of zero Lebesgue measure can be covered
by intervals of length $\ab2\ab\gamma\,$ of an arbitrarily small
total length.

\bigskip

\noindent{\bf 3. Examples and related problems.}

\m

In this section, we give only hints at the proofs; mostly the
proofs can be obtained by adapting those in Sections 1\,--\,2\ab.

\s

\noi{\bf3.0.} In possible applications to physical systems
with fast oscillating random perturbations, the
action-angle coordinates often cannot be introduced globally,
in particular, not in neighborhoods of critical points of first
integrals. In such cases our results cannot be applied immediately;
but they can after some adaptation.

Let $\bigl(\boldkey H^\eps(t),\, \boldsymbol
\Phi^\eps(t)\bigr)$ be as in Section 1; let $G$ be a region
in $\Bbb R^n$ with a smooth boundary $\partial G$. Let the
conditions of Theorem 1.1 be satisfied for $\,\boldkey h\,$
belonging to some neighborhood of the closure $\overline G$
of $G$. Let $\bigl(\hat{\boldkey H}^\eps(t),\, \hat{
\boldsymbol \Phi}^\eps(t)\bigr)$ be the stochastic process
$\bigl(\boldkey H^\eps(t),\, \boldsymbol\Phi^\eps(t)\bigr)$
stopped at the time $\,\tau^\eps\ab$ at which $\boldkey
H^\eps(t)$ leaves $G$: $\bigl(\hat{\boldkey H}^\eps(t),\, \hat{
\boldsymbol \Phi}^\eps(t)\bigr) = \bigl(\boldkey H^\eps(t
\wedge \tau^\eps),\, \boldsymbol\Phi^\eps(t \wedge \tau^\eps)
\bigr)$.

Then the process $\hat{\boldkey H}^\eps(t)$ converges in
distribution as $\,\eps \to 0\ab$ to the diffusion
process $\hat{\boldkey H}(t) $ in $\overline G$ governed by the
differential operator $L$ given by formula (1.28) -- with the boundary
condition $L\ab f(\boldkey h) = 0$ for $\,\boldkey h \in \partial G$
(this process stops at the time it leaves $G$).

The proof is a slight modification of that of Theorem 1,1, all
formulas being more complicated because we have to take integrals
truncated after the time of leaving $G$. We require the conditions
of Theorem 1.1 to be satisfied in a neighborhood of $\overline G$
rather than just in $\overline G$ because the process $\tilde{
\boldkey H}^\eps(t)$ approximating $\boldkey H^\eps(t)$ (see
 Section 1) may be outside~$\overline G $ while $\boldkey
H^\eps(t)$ is in $\overline G$.

If the conditions of Theorem 1.1 are satisfied for $\,\boldkey h\,$
not in the whole $\Bbb R^n$ but only in a
neighborhood of $\overline G$, then we can extend the coefficients
$\,\boldkey b(\boldkey h,\, \boldsymbol\varphi,\, \boldkey w)\ab$
to $\,\boldkey h\,$ outside $G$ so that the conditions of Theorem~
1.1 are satisfied in the whole space; and the process
$\hat{\boldkey H}^\eps(t)$ converges in distribution
to $\hat{\boldkey H}(t)$.

Now suppose the conditions of our Theorem 1.1 are not satisfied
in our region $G$, but there exists a sequence $G_1 \subset
G_2 \subset ... \subset G_n \subset ...$ of subregions of $G$,
such that each $G_n$ is contained in the next one with some
its neighborhood, $\bigcup_{\,n} G_n = G$, and the conditions
of Theorem 1.1 are satisfied in each $G_n$. Let $\hat{\boldkey
H}^\eps_n(t)$ be the process $\boldkey H^\eps(t)$ stopped at
the  time $\,\tau^\eps_n\ab$ of its leaving $G_n$. If the family
of function-space distributions of $\hat{\boldkey H}^\eps_n(\bullet)$
is tight, we get that there is a sequence $\,\eps_n \to 0\ab$
such that the processes $\hat{\boldkey H}^{\eps_n}_n(
t)$ converge un distribution of some stochastic
process $\hat{\boldkey H}(t) $ in $\overline G$ being a solution
of the martingale problem associated with the operator $L$;
but we don't know whether the solution of the martingale problem
is unique or if the limiting process is a Markov one.

There can be two situations here: all limiting processes starting
from an interior point of the region $G$ never reach the boundary
$\partial G$; or this boundary is accessible from interior points.
In the first case we can prove that the limiting distribution
{\it is\/} unique, and the corresponding stochastic process is
a Markov one, never leaving $G$; in the second case, under some
conditions, the processes $\hat{\boldkey H}^\eps(t)$
(stopped at reaching $\partial G$) will converge in distribution
to the process $\hat{\boldkey H}(t)$ with the generator
$L$ given by (1.28), stopping at the time $\,\tau\,$ of leaving
$G$ (but we are not interested in this now: in what we are
going to consider in Subsection 3.1 the boundary will be inaccessible).

\m

\noi{\bf3.1.} Keeping this in mind, let us consider a simple example 
of applying what was
said: a simple, but still a meaningful one. Let us start with a
non-perturbed system.

\s

A simple system describing an oscillator with one degree of
freedom is one in the two-dimensional space with solutions
moving on circles centered at $(0,\, 0)$ with its own angular
velocity on each circle. The function $H(\boldkey x) =
H(x_1,\, x_2) = \tfrac12\ab(x_1^2 + x_2 ^2)$ is a smooth first
integral of such a system; we can write the system describing
this oscillator as\linebreak  $\dot X_1(t) = -\, \omega\bigl(H(
\boldkey X(t))\bigr) \cdot X_2(t)$, $\dot X_2(t) = \omega\bigl(H(
\boldkey X(t))\bigr) \cdot X_1(t)$. If the function $\,\omega(
h)$, $\ab h \in [0,\, \infty)$, \ab is smooth up to the point $0$,
the coefficients $-\,\omega\bigl(H(\boldkey x)\bigr) \cdot x_2$,
$\omega\bigl(H(\boldkey x)\bigr) \cdot x_1\ab$ of this system
are smooth functions of $\,\boldkey x = (x_1,\, x_2)$.

\s

Now let us consider two independent oscillators; we'll have
two copies of the plane $\Bbb R^2$ with coordinates $\,\boldkey
x_1 = (x_{11},\, x_{12})\ab$ in the first plane and $\,\boldkey
x_2 = (x_{21},\, x_{22})\ab$ in the second; the functions
$\boldkey X_1(t) = \bigl(X_{!1}(t),\, X_{12}(t)\bigr)$,
$\boldkey X_2(t) = \bigl(X_{21}(t),\, X_{22}(t)\bigr)$ will be
the solutions of the system
$$
\aligned
    \dot X_{11}(t) &= -\, \omega_1\bigl(H(\boldkey X_1(t))\bigr)
    \cdot X_{12}(t),        \\
    \dot X_{12}(t) &= \te\ \omega_1\bigl(H(\boldkey X_1(t))\bigr)
    \cdot X_{11}(t),
\endaligned\qquad
\aligned
    \dot X_{21}(t) &= -\, \omega_2\bigl(H(\boldkey X_2(t))\bigr)
    \cdot X_{22}(t),        \\
    \dot X_{22}(t) &= \te\ \omega_2\bigl(H(\boldkey X_2(t))\bigr)
    \cdot X_{21}(t),
\endaligned
\eqno(3.1)
$$
where $\,\omega_1(h)$, $\omega_2(h)$, $\ab h \in [0,\, \infty)$.
\ab are smooth functions; let us assume that $\,\omega_i(h) > 0\ab$
for all $\,h \in [0,\, \infty)$, $\ab i = 1$, $2$.

In lieu of Condition $\star$ we'll have

\noi Condition $\clubsuit$: for every $\,\theta \in (0,\, \infty)\ab$
the intersection of the set $\{(h_1,\, h_2) \col T_1(h_1)/T_2(h_2)$
$=\theta\}$ with every compact subset of $(0,\, \infty)^2$ is either
empty, or consists of finitely many separate points and finitely
many smooth curves.

\s

For each of the subsystems in (3.1) we can introduce in $\Bbb R^2
\setminus \{\bold 0\}$ action-angle-type coordinates
$(h_i,\, \varphi_i)\in (0,\, \infty) \times \Bbb T^1$, related
to $\,x_{i1}$, $x_{i2}\ab$ by
$$
    x_{i1} = \sqrt{2\ab h_i} \cdot \cos(2\ab\pi\ab\varphi_i), \qquad
    x_{i2} = \sqrt{2\ab h_i} \cdot \sin(2\ab\pi\ab\varphi_i).
\eqno(3.2)
$$
In these coordinates the equations (3.1) become
$$
\aligned
    \dot H_1(t) &= 0, \\
    \dot \Phi_1(t) &= \omega_1\bigl(H_1(t)\bigr),
\endaligned\qquad\quad
\aligned
    \dot H_2(t) &= 0, \\
    \dot \Phi_2(t) &= \omega_2\bigl(H_2(t)\bigr).
\endaligned
\eqno(3.3)
$$
Here $H_i(t) = H_i\bigl(\boldkey X_i(t)\bigr)$, $\Phi_i(t)$
is the angle coordinate of $\boldkey X_i(t)$.

\s

Now let us consider fast-oscillating random perturbations of
system (3.1):
$$
\aligned
    \dot{\boldkey X}^\eps_1(t) &= \eps^{- 1} \boldkey b_1\big(
    \boldkey X^\eps_1(t),\, \boldkey X^\eps_2(t),\, w_1(t/\eps^2)\bigr), \\
    \dot{\boldkey X}^\eps_2(t) &= \eps^{- 1} \boldkey b_2\big(
    \boldkey X^\eps_1(t),\, \boldkey X^\eps_2(t),\, w_2(t/\eps^2)\bigr),
\endaligned
\eqno(3.4)
$$
where $\,w_1(s)$, $w_2(s)\ab$ are two independent Wiener processes
on the circle $\Bbb T^1$ of length~$1$, $\dsize\int_{\Bbb T^1}
\boldkey b_i(\boldkey x_1,\, \boldkey x_2,\, w_i)\ dw_i = \tilde{
\boldkey b}_i(\boldkey x_i) = \bigl(-\,\omega_i\bigl(H(\boldkey
x_i)\bigr) \cdot x_{i2},\, \omega_i\bigl(H(\boldkey
x_i)\bigr) \cdot x_{i1}\bigr)$
(for simplicity's sake, we consider the right-hand
sides~$\,\boldkey b_i\ab$ depending only on $\,w_i(t/\eps^2)$, \ab
not on the whole $\,\boldkey w(t/\eps^2) = \bigl(w_1(t/\eps^2),\,
w_2(t/\eps^2)\bigr)$).

\s

The perturbed oscillators are no longer
independent, because the right-hand sides depend not on one
$\boldkey X^\eps_i(t)$, but on both $\ab \boldkey X^\eps_1(t)$,
$\boldkey X^\eps_2(t)$.

\s

Let us introduce some notations: $\,\boldkey b_i
(\boldkey x_1,\, \boldkey x_2,\, w_i) = \bigl(b_{i1}(\boldkey x_1,\,
\boldkey x_2,\, w_i),\, b_{i2}(\boldkey x_1,\, \boldkey x_2,\, w_i)
\bigr)$, $\ab\tilde{\boldkey b}_i(\boldkey x_i)$
$= \bigl(\tilde b_{i1}(\boldkey x_i),\, \tilde b_{i2}(\boldkey x_i)\bigr)$,
$\ab\alpha_{ij}(\boldkey x_1,\, \boldkey x_2,\, w_i) = \boldkey
b_{ij}(\boldkey x_1,\, \boldkey x_2,\, w_i) - \tilde b_{ij}(\boldkey x_i)$.
Of course\linebreak $\dsize\int_{\Bbb T^1} \alpha_{ij}(\boldkey x_1,\,
\boldkey x_2,\, w_i)\ dw_i = 0$.

\s

We impose on the functions $ \,\alpha_{ij}\ab$ the following condition:

\s

\noi Condition $\clubsuit\, \clubsuit$\ab: For every $\,i = 1 $,  $2$,
and every $\,\boldkey x_1 $, $x_2\ab$ the functions $\,\alpha_{i1}(
\boldkey x_1,\, \boldkey x_2,\, w_i)$ \ab and $\,\alpha_{i2}(\boldkey x_1,\,
\boldkey x_2,\, w_i)$ \ab are linearly independent.

\s

In the action-angle coordinates the system (3.4) takes the form
$$
\aligned
    \dot H^\eps_1(t) &= \eps^{- 1} \cdot\beta_{11}\bigl(H^\eps_1(t),\,
    H^\eps_2(t),\, \Phi^\eps_1(t),\, \Phi^\eps_2(t),\,
    w_1(t/\eps^2)\bigr),            \\
    \dot\Phi^\eps_1(t) &= \eps^{- 1} \cdot\bigl[\omega_1\bigl(H^\eps_1(t)\bigr)
    + \beta_{12}\bigl(H^\eps_1(t),\,
    H^\eps_2(t),\, \Phi^\eps_1(t),\, \Phi^\eps_2(t),\,
    w_1(t/\eps^2)\bigr)\bigr],            \\
    \dot H^\eps_2(t) &= \eps^{- 1} \cdot \beta_{21}\bigl(H^\eps_1(t),\,
    H^\eps_2(t),\, \Phi^\eps_1(t),\, \Phi^\eps_2(t),\,
    w_2(t/\eps^2)\bigr),            \\
    \dot\Phi^\eps_2(t) &= \eps^{- 1} \cdot\bigl[\omega_2\bigl(H^\eps_1(t)\bigr)
    + \beta_{22}\bigl(H^\eps_1(t),\,
    H^\eps_2(t),\, \Phi^\eps_1(t),\, \Phi^\eps_2(t),\,
    w_2(t/\eps^2)\bigr)\bigr],
\endaligned
\eqno(3.5)
$$
where
$$
\aligned
    \beta_{i1}(h_1,\, h_2,\, \varphi_1,\, \varphi_2,\, w_i)
    &= \sqrt{2\ab h_i} \cdot\bigl(\cos(2\ab\pi\ab \varphi_i) \cdot \alpha_{i1}
    + \sin(2\ab\pi\ab \varphi_i) \cdot \alpha_{i2}\bigr), \\
    \beta_{i2}(h_1,\, h_2,\, \varphi_1,\, \varphi_2,\, w_i)
    &= \dfrac1{\sqrt{2\ab h_i}} \cdot\bigl(
    -\,\sin(2\ab\pi\ab \varphi_i) \cdot \alpha_{i1}
    + \cos(2\ab\pi\ab \varphi_i) \cdot \alpha_{i2}\bigr),
\endaligned
\eqno(3.6)
$$
$\,\alpha_{ij}\ab$ being taken at the point $\bigl(\sqrt{2\ab h_1}
\cdot \cos(2\ab\pi\varphi_1),\, \sqrt{2\ab h_i} \cdot \sin(
2\ab\pi\varphi_1),\, \sqrt{2\ab h_2}
\cdot \cos(2\ab\pi\varphi_2),$\linebreak $\sqrt{2\ab h_2} \cdot \sin(
2\ab\pi\varphi_2),\,w_i\bigr)$.

\s

Now we return to Section 1, to the equation $\tfrac12\,\Delta
U(\boldkey w) = -\,G(\boldkey w)$. In our present case it is
an ordinary differential equation on a circle:
$$
    \tfrac12\, U''(w) = -\, G(w),\qquad w \in \Bbb T^1.
\eqno(3.7)
$$
We still need $\dsize\int_{\Bbb T^1} G(w)\ dw = 0$, but not
H\"older continuity of $G(w)$: just continuity is enough. We
can write a solution of the equation (3.7) explicitly:
$$
    U(w) = \int_0^w (w -\tfrac12 - v)^2 \cdot G(v)\ dv
        + \int_w^1 (w + \tfrac12 - v)^2 \cdot G(v)\ dv;
\eqno(3.8)
$$
$$
    U'(w) = -\,2 \int_0^w (\tfrac12 + v) \cdot G(v)\ dv
        + 2 \int_w^1 (\tfrac12 - v) \cdot G(v)\ dv.
\eqno(3.9)
$$
The solution given by (3.8) is normalized by the condition
$\dsize\int_{\Bbb T^1} U(w)\ dw = 0$.

Let $U_{ij}(\boldkey x_1,\, \boldkey x_2,\, w_i)$ be the solution
of the equation
$$
    \dfrac12\,\dfrac{\partial^2 U_{ij}(\boldkey x_1,\,
        \boldkey x_2,\, w_i)}{\partial w_i^2}
    = -\, \alpha_{ij}(\boldkey x_1,\, \boldkey x_2,\, w_i)\ab;
\eqno(.10)
$$
then the four components of the solution $\,\boldkey u\,$ of
the equation (1.11) (in the action-angle coordinates) are
$$
\aligned
    u_{i1}(h_1,\, h_2,\, \varphi_1,\, \varphi_2,\, w_i)
    &= \sqrt{2\ab h_i} \cdot\bigl(\cos(2\ab\pi\ab \varphi_i) \cdot U_{i1}
    + \sin(2\ab\pi\ab \varphi_i) \cdot U_{i2}\bigr), \\
    u_{i2}(h_1,\, h_2,\, \varphi_1,\, \varphi_2,\, w_i)
    &= \dfrac1{\sqrt{2\ab h_i}} \cdot\bigl(
    -\,\sin(2\ab\pi\ab \varphi_i) \cdot U_{i1}
    + \cos(2\ab\pi\ab \varphi_i) \cdot U_{i2}\bigr).
\endaligned
\eqno(3.11)
$$

Now we can evaluate the coefficients $\tilde A_{ij}$, $\tilde
B_i$, $\overline A_{ij}$, $\overline B_i$ (see formulas (1.22),
(1.21), (1.26)). The integrals for the coefficients $\tilde A_{ij}$
can be taken over $\Bbb T^1$ rather than over $\Bbb T^2$:
$$
\aligned
    &\tilde A_{ii}(h_1,\, h_2,\, \varphi_1,\, \varphi_2)
    = \int_{\Bbb T^1} \bigl(\dfrac{\partial u_{i1}(h_1,\, h_2,\,
     \varphi_1,\, \varphi_2,\, w_i)}{\partial w_i}\bigr)^2\, dw_i \\
    &\qquad= \int_{\Bbb T^1} \bigl[\sqrt{2\ab h_i} \cdot
     \bigl(\cos(2\ab\pi\varphi_i)
    \cdot \dfrac{\partial U_{i1}}{\partial w_i} +
     \sin(2\ab\pi\varphi_i)\cdot \dfrac{\partial U_{i2}}
     {\partial w_i}\bigr)\bigr]^ 2\, dw_i
\endaligned
\eqno(3.12)
$$
(the functions $U_{ij}$ are taken here at the arguments $\bigl(
\sqrt{2\ab h_1}
\cdot \cos(2\ab\pi\varphi_1),\, \sqrt{2\ab h_1} \cdot \sin(
2\ab\pi\varphi_1),$\linebreak $\sqrt{2\ab h_2}
\cdot \cos(2\ab\pi\varphi_2),\, \sqrt{2\ab h_2} \cdot \sin(
2\ab\pi\varphi_2),\,w_i\bigr)$). It follows from Condition
$\clubsuit\, \clubsuit$ that the coef-\linebreak ficients
$\tilde A_{ii}$ are strictly positive for all $ \,h_i > 0$;
they go to $0$ as $\,h_i \to 0 $.

As for  $\tilde A_{12} = \tilde A_{21}$, we have:
$$
\aligned
    \tilde A_{12}(h_1,\, h_2,\, \varphi_1,\, \varphi_2)
    &= \iint_{\Bbb T ^2} \bigl[\sqrt{2\ab h_1} \cdot
     \bigl(\cos(2\ab\pi\varphi_1)
    \cdot \dfrac{\partial U_{11}}{\partial w_1} +
     \sin(2\ab\pi\varphi_1)\cdot \dfrac{\partial U_{12}}
     {\partial w_1}\bigr)\bigr] \times  \\
     &\times \bigl[\sqrt{2\ab h_2} \cdot \bigl(\cos(2\ab\pi\varphi_2)
    \cdot \dfrac{\partial U_{21}}{\partial w_2} +
     \sin(2\ab\pi\varphi_2)\cdot \dfrac{\partial U_{22}}
     {\partial w_2}\bigr)\bigr]\ dw_1\ab dw_2.
\endaligned
\eqno(3.13)
$$
This double integral is the product of two one-dimensional
integrals, both of them are equal to $0$, and $ \tilde A_{12}
= 0 $.

We see from this that Condition $\star\,\star     $ is satisfied
for $(h_1,\, h_2)$ belonging to every compact subset of $(0,\,
\infty)^2$.

The averaged matrix $\bigl(\overline A_{ij}(h_1,\, h_2)\bigr)$ is
also a diagonal one with positive diagonal entries.\linebreak When one of
$\,h_i\ab$ goes to $0$, we have:
$$
\aligned
    &\overline A_{11}(h_1,\, h_2) = D_{11}(h_2) \cdot h_1 + O(h_1^{3/2})
    \quad (h_1 \to 0), \\
    &\overline A_{22}(h_1,\, h_2) =
    D_{22}(h_1) \cdot h_2 + O(h_2^{3/2}) \quad (h_2 \to 0),
\endaligned
\eqno(3.14)
$$
where $D_{11}(h_2) > 0 $, $D_{22}(h_1) > 0$,
$$
\aligned
    &D_{11}(h_2)
    = \int_{\Bbb T^1} \Bigl[\Bigl(\dfrac{\partial U_{11}\bigl(0,\, 0,\,
    \sqrt{2\ab h_2} \cdot \cos(2\ab\pi\varphi_2),\,
    \sqrt{2\ab h_2} \cdot \sin(2\ab\pi\varphi_2),\, w_1\Bigr)}
    {\partial w_1}\bigr)^2      \\
    &\qquad\qquad\quad\ \ + \Bigl(\dfrac{\partial U_{12}\bigl(0,\, 0,\,
    \sqrt{2\ab h_2} \cdot \cos(2\ab\pi\varphi_2),\,
    \sqrt{2\ab h_2} \cdot \sin(2\ab\pi\varphi_2),\, w_1\bigr)}
    {\partial w_1}\Bigr)^2\Bigr]\ dw_1,
\endaligned
\eqno(3.15)
$$
and similarly for $D_{22}$.

\s

Formula (1.21) becomes
$$
    \tilde B_i(h_1,\, h_2,\, \varphi_1,\, \varphi_2)
    = \iint_{\Bbb T^2}\bigl(\dfrac{\partial u_{i1}}{\partial h_1}
    \cdot \beta_{11} + \dfrac{\partial u_{i1}}{\partial h_2}
    \cdot \beta_{21} + \dfrac{\partial u_{i1}}{\partial\varphi_1}
    \cdot \beta_{12} + \dfrac{\partial u_{i1}}{\partial\varphi_2}
    \cdot \beta_{22}\bigr)\ dw_1\, dw_2\ab;
\eqno(3.16)
$$
and it's easy to write an expression for $\overline B_i(h_1,\, h_2)$\ab.

\s

As for the asymptotics of $\overline B_i(h_1,\, h_2)$ as $\,h_i
\to 0$, \ab only the $\,i$-th \ab term ($i = 1$, $2$) in (3.16) is of order~$1$,
all other are $O(h_i)$\ab; and $\overline B_1(h_1,\, h_2) = \tfrac12\,
D_{11}(h_2) + O(h_1^{1/2})\ab$ \ab as $\,h_1 \to 0$, $\overline B_2(h_1,\,
h_2) = \tfrac12\, D_{22}(h_1) + O(h_2^{1/2})$ \ab as $\,h_2 \to 0$ (the
fact that $\botsmash{\dsize\int_{\Bbb T^1} U_{ij} \cdot \alpha_{ij}\ dw_i}
= \dsize\int_{\Bbb T^1} \bigl(\dfrac{\partial U_{ij}}{\partial w_i}
\bigr)^2\, dw_i$ \,is used).

\s

It follows from this that the boundary of the region $G = (0,\, \infty)^2$
is inaccessible for the limiting process starting from interior
points of $G$:

\s

{\bf Lemma 3.1.} {\it Let the stochastic process\/ $\bigl(H_1(t),\,
H_2(t)\bigr)$ {\rm(not necessarily a diffusion one)} be a
solution of the martingale problem associated
with the differential operator
$$
    Lf(h_1,\, h_2)
    = \dfrac12 \sum_i A_{ii}(h_1,\, h_2) \cdot
        \dfrac{\partial^2 f}{\partial h_i^2}
    + \sum_i B_i(h_1,\, h_2) \cdot \dfrac{\partial f}{\partial h_i};
\eqno(3.17)
$$
and\/ $B_i(h_1,\, h_2) > A_{ii}(h_1,\, h_2)/2\ab h_i - O(h_i^{1/2})$
as\/ $\,h_i \to 0.$ \ab Then the line\/ $\{(h_1,\, h_2)\col
h_i = 0\}$ is inaccessible for this process\ab}.

\s

The {\bf proof} is based on the fact that for the function
$\,f(h_i) = \ln(-\, \ln h_i)\ab$ we have $(Lf)(h_1,\, h_2)
< 0$ for sufficiently small $\,h_i$.

\s

By what was said in Subsection 3.0 (with $(0,\, \infty)^2$
as $G$ and $(\delta,\, \infty)^2$ with the ``rounded'' corner
as $G_n$), the process $\boldkey H^\eps(t)$ converges
in distribution to the diffusion process
in $(0,\, \infty)^2$ with the generating operator~$L$.

\m

The systems (3.1) can be rewritten as Hamiltonian systems
with energy functions $H_1(\boldkey x)$, $H_2(\boldkey x)$ depending
only on $|\boldkey x|$ and having only one critical point each.

\s

If the first integrals of the two oscillators are not $\tfrac12\,
(x_1^2 + x_2^2)$ but some smooth functions $H_1$, $H_2$ with
one critical point each at which the matrix of second derivatives
is positive definite, and $\,\lim_{|\boldkey x| \to \infty}
H_i(\boldkey x) = \infty$, \ab we get similar results with the
quarter-plane $\bigl(\min_{\ab\boldkey x} H_1(\boldkey x),\, \infty
\bigr) \times \bigl(\min_{\ab\boldkey x} H_2(\boldkey x),\, \infty
\bigr)$ replacing $(0,\, \infty)^2$, only we can no longer
write so explicit formulas as we did for $H(\boldkey x) = \tfrac12
\,(x_1^2 + x_2^2)$\ab.

\m

\noi{\bf 3.2.} Consider a system of two independent
one-degree-of-freedom oscillators:
$$
    \dot{\boldkey X}_1(t)
    = \overline\nabla H_1\bigl(\boldkey X_1(t)\bigr), \qquad
    \dot{\boldkey X}_2(t)
    = \overline\nabla H_2\bigl(\boldkey X_2(t)\bigr).
\eqno(3.18)
$$
We assume that the Hamiltonians $H_1$ and $H_2$ are smooth
enough, $\ab\lim_{\ab |\boldkey x| \to \infty} H_i(\boldkey
x)$\linebreak $= \infty$, $\ab i = 1$, $2$, \ab and each of these functions
has a finite number of non-degenerate critical points. To be
specific, let $H_1(\boldkey x)$ have just one minimum at
$\ab O \in \Bbb R^2$, $\nabla H_1(\boldkey x) \neq \bold 0$ for
$\,\boldkey x \neq O$; and $H_2(\boldkey x)$ three critical
points: minima at $O_1$ and $O_2$, and a saddle point at $O_3$,
$\nabla H_2(\boldkey x) \neq \bold0$ for $\,\boldkey x \neq
O_1$, $O_2$, $O_3$. 

For the first equation in (3.18) we can introduce the action-angle
coordinates in the whole plane minus one point $O$.
All solutions $\boldkey X_1(t)$ of the first equation (3.18)
are periodic functions moving on a level set $\{\boldkey
x \in \Bbb R^2\col H_1(\boldkey x) = h\}$. The set of all
level sets can be parametrized by the graph~$\Gamma_1$ consisting
of one vertex $\Cal O$ corresponding to $O$ and one edge~$I$
with coordinate $\,h \in \bigl[H_1(O),\, \infty\bigr)$. Define
the mapping $Y_1: \Bbb R^2 \mapsto \Gamma_1$ taking as $Y_1(
\boldkey x)$ the point of the graph with coordinate $\,h =
H_1(\boldkey x)$.

\s

The Hamiltonian $H_2$ has two wells. The level set $\{\boldkey x\col
H_2(\boldkey x) = H_2(O_3)\}$ is figure-of-eight shaped. It divides
the plane into three open regions: $G_1$, inside one loop of the
figure eight, containing the
equilibrium point $O_1$; $G_2$, inside the other loop containing
$O_2\ab$; and $G_3$, outside the figure of eight. In each of the
regions $G_1 \setminus \{O_1\}$, $G_2 \setminus \{O_2\}$, and
$G_3$ we can introduce its own action-angle coordinates. Identifying
points of each connected component of the level sets $\{\boldkey x
\col H_2(\boldkey x) = h\}$, we get a graph $\Gamma_2$ with three
vertices $\Cal O_1$, $\Cal O_2$, $\Cal O_3$ corresponding to
$O_1$, $O_2$ and the level set $\{\boldkey x\col H_2(\boldkey x)
= H_2(O_3)\}$, and three edges: $I_1$ corresponding to the
connected components of level sets lying in $G_1$, $I_2$ and
$I_3$ the same with $G_2$, $G_3$. The edges $I_1$, $I_2$, $I_3$
are glued together at the vertex $\Cal O_3$. We can introduce
coordinates on the graph  $\Gamma_2$: if a point $\, y\,$ of
the graph belongs to $I_k$ and it corresponds to a connected
component of $\{\boldkey x \col H_2(\boldkey x) = h\}$, we take
as its coordinates $(h,\, k)$ ($h\,$ taken alone cannot serve
as the coordinate for graph points, because there are points
on the edges $I_1$ and $I_2$ with the same $\,h$). The second
coordinate can be considered as a discrete-valued first integral
of our system.

Define a mapping $Y_2: \Bbb R^2 \mapsto \Gamma_2$\ab: $Y_2(
\boldkey x) = (h,\, k)$ if $H_2(\boldkey x) = h$, and $\,
\boldkey x$ (and the whole solution $\boldkey X_2(t)$ starting
at $\,\boldkey x$) belongs to $G_k$.

Let $\ab\Pi = \Gamma_1 \times \Gamma_2$\ab; a geometric object
of this kind is called ``an open book'' (see [10], Appendix).
The open book
$\Pi$ has three pages $\,\pi_k = \Gamma_1 \times I_k$, $\ab
k = 1$, $2 $, $3$, \ab and the binding $\Gamma_1 \times\{
\Cal O_3\}$. We assume that Vondition $\star$ is satisfied 
on each page.

Construction of the graph corresponding to a one-degree-of-freedom
Hamiltonian and of the open book in the case of many conservation
laws is described in detail in [7], Chapters 8, 9.

Define the mapping $\boldkey Y\!: \Bbb R^2 \times \Bbb R^2 \mapsto
\Pi$\ab: $\boldkey Y(\boldkey x_1,\, \boldkey x_2) = \bigl(Y_1(
\boldkey x_1),\, Y_2(\boldkey x_2)\bigr)$.

\s

Now let us consider fast oscillating perturbations of the system
(3.18):
$$
\aligned
    \dot{\boldsymbol\xi}^\eps_1(t)
    &= \overline\nabla H_1\bigl(\boldsymbol\xi^\eps_1(t)\bigr)
    + \boldsymbol\alpha_1\bigl(\boldsymbol\xi^\eps_1(t),\,
    \boldsymbol\xi^\eps_2(t),\, w_1(t/\eps)\bigr),  \\
    \dot{\boldsymbol\xi}^\eps_2(t)
    &= \overline\nabla H_2\bigl(\boldsymbol\xi^\eps_2(t)\bigr)
    + \boldsymbol\alpha_2\bigl(\boldsymbol\xi^\eps_1(t),\,
    \boldsymbol\xi^\eps_2(t),\, w_2(t/\eps)\bigr).
\endaligned
\eqno(3.19)
$$
Here $\,\boldsymbol \alpha_1$,
$\boldsymbol \alpha_2\ab$ are smooth functions with $\dsize
\int_{\Bbb T^1} \boldsymbol\alpha_i(\boldkey x_1,\, \boldkey x_2,
\, w)\ dw = \bold 0\ab$ for all $\,\boldkey x_1$, $\boldkey x_2
\in \Bbb R^2$, $\ab i = 1$, $2$. Then we can check that $\bigl(
\boldsymbol\xi^\eps_1(t),\, \boldsymbol\xi^\eps_2(t)\bigr)$
converges in probability as $\,\eps \to 0$, \ab uniformly on
every finite interval $[0,\, T]$, to the trajectory $\bigl(
\boldkey X_1(t),\, \boldkey X_2(t)\bigr)$ of the non-perturbed
system (3.18) with the same initial condition; so we can
consider $\bigl(\boldsymbol\xi^\eps_1(t),\,
\boldsymbol\xi^\eps_2(t)\bigr)$ as the result of small perturbations
applied to the system (3.18). Since $\ab H_1(\boldkey x_1)$,
$H_2(\boldkey x_2)\ab$ are first integrals of (3.18), this
implies that $\bigl(H_1\bigl(\boldsymbol\xi^\eps_1(t)\bigr),\,
H_2\bigl(\boldsymbol\xi^\eps_2(t)\bigr)\bigr)$ converges to
the constant $\bigl(H_1(\boldkey x_1),
\, H_2(\boldkey x_2)\bigr)$ as $\,\eps \to 0\ab$ for every
$\,t \in [0,\, \infty)$.

\s

But on large time intervals of order of $\,\eps^{- 1}$ the
deviation of $H_i\bigl(\boldsymbol\xi^\eps_i(t)\bigr)$ from
its initial value can be of order $1$. We take $\boldkey
X_i^\eps(t) = \boldsymbol\xi_i^\eps(t/\eps)$. Assume, for simplicity,
that the vector fields $\,\boldsymbol \alpha_i$ are also Hamiltonian:
$\ab\boldsymbol\alpha_i(\boldkey x_1,\, \boldkey x_2,\, w)
= \overline \nabla_{\boldkey x_i} \Cal H_i(\boldkey x_1,\,
\boldkey x_2,\, w)$, where the functions $\Cal H_i$, $\ab i =
1$, $2$, \ab are smooth enough, $\botsmash{\dsize\int_{\Bbb T^1}}
\Cal H_i(\boldkey x_1,\, \boldkey x_2,\, w)\ dw =0$. Then
$\boldkey X_1^\eps(t)$, $\boldkey X_2 ^\eps(t)$ satisfy
the equations
$$
\aligned
    \dot{\boldkey X}_1^\eps(t) &
    = \eps^{- 1} \overline\nabla_{\boldkey x_1} \bigl[
    H_1\bigl(\boldkey X_1^\eps(t)\bigr)
    + \Cal H_1\bigl(\boldkey X_1^\eps(t),\, \boldkey X_2^\eps(t),
    w_1(t/\eps^2)\bigr)\bigr],      \\
    \dot{\boldkey X}_2^\eps(t) &
    = \eps^{- 1} \overline\nabla_{\boldkey x_2} \bigl[
    H_2\bigl(\boldkey X_2^\eps(t)\bigr)
    + \Cal H_2\bigl(\boldkey X_1^\eps(t),\, \boldkey X_2^\eps(t),
    w_2(t/\eps^2)\bigr)\bigr].
\endaligned
\eqno(3.20)
$$

Let $\boldkey Y^\eps(t) = \bigl(Y_1(t),\, Y_2(t)\bigr) =
\boldkey Y\bigl(\boldkey X^\eps_1(t),\,
\boldkey X^\eps_2(t)\bigr)$ be the projection of $\bigl(\boldkey
X^\eps_1(t),\, \boldkey X^\eps_2(t)\bigr)$ onto the open book~
$\Pi$; let $\,\tau^\eps\ab$ be the first time that the process
$\boldkey Y^\eps(t)$ hits the boundary of the page where it
started. Let $\hat{\boldkey Y}^\eps(t)$ be the process $\boldkey Y^\eps(t)$
stopped at the time $\,\tau^\eps$: $\ab\hat{\boldkey Y}^\eps(t)
= \boldkey Y^\eps(t \wedge \tau^\eps)$. Keeping in mind Theorem 1.1
and what was said in Subsections 3.0, 3.1, we can expect that
the time $\,\tau^\eps\ab$ is, in fact. the time of reaching the
binding (the parts $\{\Cal O\} \times \Gamma_2$, $\Gamma_1 \times
\{\Cal O_1\}$, $\Gamma_1 \times \{\Cal O_2\}$ of the boundary
being inaccessible) and that under
some natural additional assumptions, in particular,
one concerning the smallness of the resonance set, the process
$\hat{\boldkey Y}^\eps(t)$ converges in distribution to a
diffusion process on the page $\,\pi_k\ab$
in which the process started (also stopping at reaching the
binding). We can hope that the {\it un\,}stopped process
$\boldkey Y^\eps(t)$ converges in distribution to some diffusion process
$\boldkey Y(t) = \bigl(Y_1(t),\, Y_2(t)\bigr)$ on the
whole open book~$\Pi $ (for white-noise-type
perturbations such results an be found in [7], Chapters
8, 9); but to which of diffusion processes on that space?

To identify the limiting diffusion process on the open book
we have to describe its behavior after hitting the binding.
One does this by calculating the generator $A$ of this process,
including its domain of definition.

Inside each page $\,\pi_k$ \ab the operator $A$ is expressed
on smooth functions $\, u(h_1,\, h_2,\, k)\ab$ belonging to
its domain of definition by the same differential operator
$L_k$ as for the process that stops at leaving $\,\pi_k$\ab, \ab
and it can be found using the approach of Subsections 3.0 and
3.1. To find the gluing conditions at the binding of $\Pi$
defining the domain of $A$ (we don't meed any boundary conditions
at the inaccessible parts of the boundary), one can use the
following observation: Consider the $6$-dimensional diffusion
process $\boldkey Z^\eps (t) =$
$\bigl(\boldkey X_1^\eps(t),\, \boldkey X_2^\eps(t),\,
w_1(t/\eps^2),\, w_2(t/\eps^2)\bigr)$ in $\ab\Cal E = \Bbb R^2 \times
\Bbb R^2 \times \Bbb T^2$.  One can write down the forward
Kolmogoroc equation for this process and see that the Lebesgue
measure $\,\lambda_6\ab$ on $\Cal E$ is invariant for $\boldkey
Z^\eps(t)$for every $\,\eps > 0$ (here we are using the Hamiltonian
form of the perturbations). The process $\boldkey Y^\eps(t)$ is
the projection of $\boldkey Z^\eps(t)$ onto $\Pi$; this implies
that the projection $\,\mu\,$ of the Lebesgue measure $\,\lambda_6\ab$
onto $\Pi$ (and we can write the density of the measure $\,\mu\,$
explicitly)is an invariant measure for the limiting process
$\boldkey Y(t)$ (assuming that it is a continuous Markov process).
This allows us to find the gluing conditions.

Consider first the process $\boldkey Y(t)$ with its first component
``frozen'': $Y_1(t) = h_1) \equiv h_1$. If the second component
$Y_2^{h_1}(t)$ of this modified process is a continuous Markov
process on $\Gamma_2$. we can calculate the generator of this
process within the edges $I_k$ of the graph $\Gamma_2$. All
possible gluing conditions at the vertex $\Cal O_3 \in \Gamma_2$
for a Markov process with continuous trajectories are described
in [7], Chapter 8. Taking into account the fact that we
know the invariant measure of the process, we can write the gluing
conditions as
$$
    \sum_{k = 1}^2 \gamma_k(h_1) \cdot
    \dfrac{\partial u\bigl(h_1,\, H_2(O_3)^-,\, k\bigr)}{\partial h_2}
    = \gamma_3(h_1) \cdot \dfrac{\partial u\bigl(h_1,\, H_2(O_3)^+,\, 3\bigr)}
    {\partial h_2}
\eqno(3.21)
$$
(remember, we are keeping $\,h_1\ab$ constant),
where $\,\gamma_k(h_1) > 0$, $\ab k = 1$, $2$, $3$, $\ab \gamma_1
(h_1) + \gamma_2(h_1) = \gamma_3(h_1)$; the coefficients
$\,\gamma_j(h_1)\ab$ can be written explicitly through the
density of the invariant measure near the point $(h_1,\,
\Cal O_3)$ and the perturbations.

Now, if we ``unfreeze'' $Y_1(t) $, the gluing conditions will
have the same form, but they are to be satisfied not for one
fixed  $\,h_1 $, \ab but for all $\,h_1 \geq H_1(O)$. This can ve
be derived from the fact that the process started near the
binding hits it in a very short time, and its first component
will not deviate much from its initial value in this short time.

Of course these arguments are far from a rigorous proof.
Much remains to be checked about the pre-limit process. In
particular, very essential for the possible limit being a
Markov process is the following property of the pre-limit
process: Let $D$ be a neighborhood of the binding of our
open book, $\ab\tau^\eps_D = \min\{t\col \boldkey Y^\eps(t)
\notin D\}$; let the initial point $\boldkey Y^\eps(0) = 
 (h_1,\, h_2,\, k)$. Then the distribution of $\boldkey T^\eps(
\tau^eps_D)$ between the pages $\,\pi_1$, $\pi_2$, $\pi_3\ab$
depends only very little on the number $\, k\,$ of the page 
it started from if $|h_2 - H_2(O_3)|$ and~$\,\eps\,$ are 
small enough (compare with Lemma 8.3.6 of [7]).

\m

\noi{\bf 3.3.} Consider now perturbations of the Landau --
Lifshitz magnetization equation
$$
    \dot{\boldkey X}(t)
    = \boldkey X(t) \times \nabla \tilde G\bigl(\boldkey X(t)\bigr).
\eqno(3.22)
$$
Here $\tilde G(\boldkey x)$, $\boldkey x \in \Bbb R^3$, is a smooth
function. It's easy to check that the flow $\boldkey X(t)$
in $\Bbb R^3$ defined by (3,22) has two first integrals
$M(\boldkey x) = |\boldkey x|^2/2\ab$ and $\tilde G(\boldkey x)$,
$\text{div}\,[\nabla M(\boldkey x) \times \nabla \tilde G(\boldkey x)] = 0$,
so that $\boldkey X(t)$ preserves the volume in $\Bbb R^3$.

From the physical perspective, it is natural to consider perturbations
of equation (3.22) that preserve the first integral $M$
(see [2]). Additive white-noise-type perturbations
were considered in [2], [6]. The results
in this case are very similar to perturbations of Hamiltonian
systems with one degree of freedom [7].

Consider now fast oscillating perturbations preserving $M(\boldkey x)$:
$$
    \dot{\boldsymbol \xi}^\eps(t)
    = \boldsymbol\xi^\eps(t) \times \nabla_{\boldkey x} G\bigl(
        \boldsymbol\xi^\eps(t),\, w(t/\eps)\bigr),
\eqno(3.23)
$$
where, as above, $\ab w(t)\ab$ is the Wiener process on the
unit-length circle $\Bbb T^1$. We assume that $\dsize\int_{\Bbb T^1}
G(\boldkey x,\, w)\ dw \equiv \tilde G(\boldkey x)$. Then one can derive
from the standard averaging principle that $\,\boldsymbol\xi^\eps(t)\ab$
converges weakly on each finite time interval $[0,\, T]$ to the
solution of (3.22) with the same initial condition. Since
$\tilde G(\boldkey x)$ is a first integral of (3.22), this implies
that $\tilde G\bigl(\boldsymbol\xi^\eps(t)\bigr)$ converges as $\,\eps \to 0\ab$
to a constant, namely to the value of the function $\tilde G$ at the
initial point. But it can be derived from Theorem 1.1 that
on time intervals of order of $\,\eps^{- 1}\ab$ it deviates from
this constant value by a distance of order $1$.

More precisely: Let $\boldkey X^\eps(t) = \boldsymbol\xi^\eps(t/\eps)$.
Let $\tilde G(\boldkey x)$ have on the sphere $S_{z_1} =$\linebreak $\{\boldkey x \in
\Bbb R^3\col M(\boldkey x) = z_1\}$ one minimum at the point
$\,\boldkey x_1 \in S_{z_1}\ab$ and one maximum at the point
$\,\boldkey x_2 \in S_{z_1}$. In order for everything to be similar
to what was considered in Sections 1\,--\,2 and Subsections 3.0\,--\,
3.1, we need to impose some conditions on the unperturbed system,
that is, on the function $\tilde G$l and on the perturbation function
$\,\boldsymbol\alpha(\boldkey x,\, w) = \nabla_{\boldkey c}
G(\boldkey x,\, w) - \nabla \tilde G(\boldkey x)$.

As the condition on the unperturbed system we take the condition
$\nabla M(\boldkey x) \times \nabla \tilde G(\boldkey x)$ $\neq \bold 0$
for any $\,\boldkey x \in S_{z_1}$, $\ab\boldkey x \neq \boldkey x_1$,
$\boldkey x_2$\ab. Under some conditions on $\,\boldsymbol
\alpha(\boldkey x,\, w)\ab$ similar to those in Section 1 and
Subsection 3.1 the process $\tilde G\bigl(\boldkey X^\eps(t)\big)$
converges in distribution to some diffusion process on the interval $\bigl(
\tilde G(\boldkey x_1),\, \tilde G(\boldkey x_2)\bigr)$. The ends of this
interval are inaccessible for the limiting process; its generator
can be calculated using Theorem 1.1 (taking into account
what was said in Subsection 3.0).

If the function $\tilde G(\boldkey x)$ on $S_{z_1}$ has some
local maxima and minima, and some saddle points,
one can obtain a limiting process on the associated graph, and
at its ``interior'' vertices some gluing conditions should be prescribed.

\s

We don't need any assumptions of smallness of the resonance type
in the case of perturbations of one equation (3.22). Consider
now a system of two equations
$$
\aligned
    \dot{\boldkey X}^\eps_1(t)
    &= \eps^{- 1} \cdot \boldkey X^\eps_1(t) \times
    \nabla_{\boldkey x_1} G_1\bigl(\boldkey X^\eps_1(t),\,
    \boldkey X^\eps_2(t), \boldkey w(t/\eps^2)\bigr),       \\
    \dot{\boldkey X}^\eps_2(t)
    &= \eps^{- 1} \cdot \boldkey X^\eps_1(t) \times
    \nabla_{\boldkey x_2} G_2\bigl(\boldkey X^\eps_1(t),\,
    \boldkey X^\eps_2(t), \boldkey w(t/\eps^2)\bigr),
\endaligned
\eqno3.24)
$$
where $\,\boldkey w(t)\ab$ is the Wiener process on the torus
$\Bbb T^m$; assume that
$$
    \int_{\Bbb T^m} G_i(\boldkey x_1,\, \boldkey x_2,\, \boldkey w)\ d\boldkey w
    = \tilde G_i(\boldkey x_i). \qquad i = 1, 2.
\eqno(3.25)
$$
Let $\Gamma_1$ and $\Gamma_2$ be the graphs counting connected
components of the level sets of the functions $\tilde G_1$,
$\tilde G_2$ on the spheres $|\boldkey x_1|^2/2 = z_1$,
$|\boldkey x_2|^2/2 = z_2$ respectively, $\Pi = \Gamma_1 \times
\Gamma_2$. Then, under the assumptions similar to those of
Theorem 1.1 (including the assumption concerning the smallness
of the resonance set), one can expect that the projection
of $\bigl(\boldkey X^\eps_1(t),\, \boldkey X^\eps_2(t)\bigr)$
onto $\Pi$ converges in distribution to a diffusion process
on $\Pi$.

\m

\noi{\bf 3.4.} The approach used in Theorem 1.1 can be applied
in a more general case when the fast-oscillating noise depends
on $\boldkey X^\eps(t)$\ab:
$$
    \dot{\boldkey X}^\eps(t)
    = \eps^{- 1} \cdot \boldkey b_1\bigl(\boldkey X^\eps(t),\,
    {\boldsymbol\zeta}^\eps(t)\bigr), \qquad
    \dot{\boldsymbol\zeta}^\eps(t)
    = \eps^{- 2} \cdot \boldkey b_2\bigl(\boldkey X^\eps(t)\bigr)
    + \eps^{- 1} \cdot \sigma\bigl(\boldkey X^\eps(t)\bigr) \,
    \dot{\boldkey w}(t).
\eqno(3.26)
$$

\bigskip

\noindent{\bf References.}

\m

1. V.I.\,Arnold. {\sl Mathematical Methods of Classical Mechanics},
Springer, 1978.

\s

2. G.\, Bertotti, I.\, Mayergoz, C.\, Serpico. {\sl Nonlinear
Magnetization Dynamics in Nano\-systems}, Elsevier, 2009.

\s

3. A.N.\,Borodin. Limit theorems for solutions of differential
equations with random right-hand side, {\sl Theory of Probability
and Applications}, {\bf 23}, 3, 1977, pp.\, 482\,--\,497.

\s

4. A.N.\,Borodin, M.I.\,Freidlin. Fast oscillating random
perturbations of dynamival systems with conservation laws,
{\sl Ann.\, Inst.\, Henri Poincar\'e}, {\bf 31}, 3, 1995, pp.\,
485\,--\,520.

\s

5. R.Cogburn, J.A.\, Ellison. Stochastic theory of adiabatic
invariance, {\sl Comm.\, in Math.\linebreak Phys.}, 149, 1992, pp.\,
97\,--\,126.

\s

6. M.\, Freidlin, W.\, Hu. On perturbations of generalized
Landau-Lifshitz dynamics. {\sl Jourmal of Stat.\, Phys.}, {\bf 144}, 5,
2011, pp.\, 978\,--\,1008.

\s

7. M.I.\,Freidlin, A.D.\ab Wentzell. {\sl Random Perturbations
of Dynamical Systems}, Sprin\-ger, 2012.

\s

8. I.I.\, Gikhman, A.V.\, Skorohod. {\sl Introduction to the
Theory of Random Processes}, W.B.\, Sanders, Philadelphia,
1969.

\s

9. R.Z.\,Khasminskii. A limit theorem for solutions of differential
equations with random right-hand side, {\sl Theory of Probability
and Applications}, {\bf 11}, 3, 1966, pp.\, 390\,--\,406.

\s

10. A.\, Ranicki. {\sl High Dimensional Knot Theory}, Springer, 1998.

\s

11. R.L.\, Stratonovich. Conditional Markov Processes and
Their Applications in the Theory of Optimal Control, {\sl Modern
Analytic and Computational Methods in Science and Mathematics},
7, American Elsevier: New York, 1968.

\end